\documentclass[12pt,a4paper,reqno]{amsart}
\usepackage[T1]{fontenc}
\usepackage[utf8]{inputenc}
\usepackage{lmodern,microtype}
\usepackage[margin=1.2in,heightrounded]{geometry}
\usepackage{amsmath,amssymb,mathtools,enumitem}
\usepackage{xcolor,needspace}
\usepackage[colorlinks=true,linkcolor=blue!55!black,citecolor=blue!55!black,urlcolor=blue!55!black]{hyperref}
\hypersetup{pdftitle={Periodic orbit asymptotics and dynamical zeta functions for smooth surface diffeomorphisms},
 pdfauthor={Gang Liao and Yao Tong},
 pdfsubject={Smooth surface dynamics, periodic orbits, Lyapunov exponents}}
\setlist[enumerate]{leftmargin=*,itemsep=3pt,topsep=4pt}
\allowdisplaybreaks
\newtheorem{theorem}{Theorem}[section]
\newtheorem{proposition}[theorem]{Proposition}
\newtheorem{lemma}[theorem]{Lemma}
\newtheorem{corollary}[theorem]{Corollary}
\theoremstyle{remark}
\newtheorem*{remark}{Remark}
\newcommand{\R}{\mathbb R}
\newcommand{\PP}{\mathcal P}
\newcommand{\F}{\mathcal F}
\newcommand{\HH}{\mathcal H}
\newcommand{\wh}{\widehat}
\newcommand{\one}{\mathbf 1}
\newcommand{\leqs}{\leqslant}
\newcommand{\geqs}{\geqslant}
\DeclareMathOperator{\Per}{Per}
\DeclareMathOperator{\Fix}{Fix}
\DeclareMathOperator{\per}{per}
\newcommand{\htop}{h_{\mathrm{top}}}
\title[Periodic orbit asymptotics]{Periodic orbit asymptotics and dynamical zeta functions for smooth surface diffeomorphisms}

\author{Gang Liao\textsuperscript{1}}
\author{Yao Tong\textsuperscript{2}}

\thanks{\raggedright\hspace*{-\normalparindent}%
\textsuperscript{1}\,School of Mathematical Sciences,
Center for Dynamical Systems and Differential Equations,
Soochow University, Suzhou 215006, China.}
\thanks{\raggedright\textsuperscript{2}\,School of Mathematical Sciences,
Peking University, Beijing 100871, China.}

\thanks{\raggedright E-mail addresses:
\textsuperscript{1}\,\href{mailto:lg@suda.edu.cn}{\nolinkurl{lg@suda.edu.cn}},
\textsuperscript{2}\,\href{mailto:2401110024@stu.pku.edu.cn}{\nolinkurl{2401110024@stu.pku.edu.cn}}.}

\date{}

\begin{document}
\begin{abstract}
Let \(f\) be a \(C^\infty\) diffeomorphism of a closed surface $M$,
and let \(\HH\) be a closed homoclinic class with topological entropy
\(h=h_{\mathrm{top}}(f|\mathcal{H})>0\) and homoclinic period \(\ell\). We prove that
\[
 \forall\,0<\chi<h,\qquad
 \lim_{\substack{n\to\infty\\\ell\mid n}}
 \frac{\#\Per_\chi(n,\HH)}{e^{nh}}=\ell,
\]
where \(\Per_\chi(n,\HH)\) denotes the saddle points in
\(\HH\) of least period \(n\) whose Lyapunov exponents lie
outside \([-\chi,\chi]\).
We also prove exponential remainders for the class and global
counting formulas. 
As a result, for every \(0<\chi<\htop(f)\), the Artin--Mazur zeta
function of all \(\chi\)-hyperbolic periodic orbits, including
sinks and sources, extends meromorphically beyond its disk of
convergence. Its leading poles are determined by the periods
of the ergodic maximal measures. This answers the meromorphic
continuation question of Buzzi \cite{BuzziICM}.
\end{abstract}
\maketitle
\hypersetup{linkcolor=black}
\tableofcontents
\hypersetup{linkcolor=blue!55!black}

\section{Introduction}

The prime orbit problem asks for precise asymptotics for
periodic orbits counted by their least period. Its classical
origins lie in the work of Selberg \cite{Selberg} and Huber
\cite{Huber} on closed geodesics on compact surfaces of constant
negative curvature. Margulis extended this theory to compact
manifolds of negative curvature \cite{Margulis1969,Margulis},
and Parry and Pollicott \cite{ParryPollicott} established the
prime orbit theorem for topologically weakly mixing basic sets
of Axiom A flows. In discrete time, the hyperbolic framework
of Anosov and Smale, together with the construction of Markov
partitions by Sinai and Bowen, provides the corresponding
finite-state symbolic model
\cite{Anosov,BowenMarkov,Bowen,Sinai,Smale}.
In particular, for a hyperbolic basic set of entropy \(h>0\)
and period \(\ell\), the number of periodic points of least
period \(n\) is asymptotic to \(\ell e^{nh}\) along multiples
of \(\ell\).

Beyond uniform hyperbolicity, logarithmic growth and precise asymptotics represent distinct aspects of dynamical complexity. For a \(C^{1+\alpha}\)
diffeomorphism \(f\) of a closed surface, with \(\alpha>0\) and
\(h=\htop(f)>0\), Katok's horseshoe theorem \cite{Katok},
based on Pesin theory \cite{Pesin}, gives
\[
     h\leqs\limsup_{n\to\infty}\frac1n\log\#\Fix(f^n),
 \qquad \Fix(f^n)=\{x:f^nx=x\}.
\]
Katok later asked whether the stronger bound
\begin{equation}\label{katok}
 \limsup_{n\to\infty}e^{-nh}\#\Fix(f^n)>0
\end{equation}
holds in the \(C^{1+\alpha}\) or \(C^\infty\) category
\cite[Problem~7.4]{KatokConjecture}.
Using countable-state symbolic
coding, Sarig \cite{Sarig} proved
a positive normalized lower bound in \eqref{katok} for hyperbolic fixed points
of \(f^n\) along multiples of a suitable integer whenever a
positive-entropy measure of maximal entropy exists. In the \(C^\infty\) category, Newhouse's entropy upper
semicontinuity theorem \cite{Newhouse} ensures the existence
of a measure of maximal entropy. Together, these results resolve the smooth
surface case of Katok's question.

Finite regularity does not suffice to deduce a positive lower bound in \eqref{katok}. Ding and Tong \cite{DT}
constructed, for every \(1<r<\infty\) and every nonnegative
sequence \(a_n=o(n)\), a \(C^r\) sphere diffeomorphism with
\(h=\htop(f)>0\) such that
\[
 \lim_{n\to\infty}\frac1n\log\#\Fix(f^n)=h,
 \qquad
 \sum_{n\geqs1}e^{-nh+a_n}\#\Fix(f^n)<\infty.
\]
The examples have a mixing homoclinic class of entropy \(h\)
containing all but finitely many periodic orbits, and every
ergodic invariant measure is hyperbolic. Thus even on such a
class, normalized counts may tend to zero after any prescribed
subexponential correction. Katok's logarithmic lower bound is
therefore optimal in finite regularity when no measure of
maximal entropy is assumed.

In general, upper bounds for the normalized counts in
\eqref{katok} require a further restriction: even the number
of saddle periodic points of a \(C^\infty\) surface
diffeomorphism can grow superexponentially
\cite{Kaloshin,Turaev,Berger}.
We therefore restrict the count to the \(\chi\)-hyperbolic saddle points.

From now on, \(M\) is a closed smooth Riemannian
surface, \(d_M\) is its Riemannian distance, and
\(f: M\to M \) is a $C^{\infty}$ diffeomorphism.
Let \(\per(x)\) denote the least period of a periodic point \(x\) of $f$.
A periodic point \(x\) is a saddle if \(D_xf^{\per(x)}\) has
one eigenvalue inside and one outside the unit circle.
Let \(E^s(x),E^u(x)\) be their respective eigenspaces, and put
\[
 \lambda^\sigma(x)=\frac1{\per(x)}
 \log\bigl\|D_xf^{\per(x)}|_{E^\sigma(x)}\bigr\|,
 \qquad \sigma\in\{s,u\}.
\]
An ergodic invariant probability measure \(\mu\) is
\emph{hyperbolic of saddle type} if its Lyapunov exponents satisfy
\(\lambda^s(\mu)<0<\lambda^u(\mu)\).
We call a saddle periodic point \(\chi\)-hyperbolic if
\(\lambda^u(x)>\chi\) and \(\lambda^s(x)<-\chi\).
For an invariant set
\(K\subset M\) and \(\chi>0\), put
\[
 \Per_\chi(n,K)
 =\{x\in K:\per(x)=n,\ \lambda^u(x)>\chi,\ \lambda^s(x)<-\chi\},
\]
where only saddle points are included. Each orbit contributes
exactly \(n\) points to this set.

A measure of maximal entropy, or a maximal measure, on a compact
invariant set \(K\) is an invariant Borel probability measure \(\mu\)
supported on \(K\) such that its metric entropy attains the topological entropy: \(h_\mu(f)=\htop(f|K)\).

For every \(0<\chi<\htop(f)\), Buzzi
\cite[Theorem~1.2]{Buzzi} refined Sarig's coding to an injective
coding and obtained the sharp lower coefficient \(p\) along
multiples of a suitable integer \(p\).
Burguet \cite{Burguet} proved the matching upper bound at
the exponential scale. Together, these results show that the
logarithmic upper growth rate of \(\#\Per_\chi(n,M)\) equals
\(\htop(f)\). Burguet's local periodic growth estimate also
identifies weak-\(*\) limits of counting measures as measures of maximal entropy
when the periodic families attain this growth rate.
These results leave open whether the normalized counts
\(e^{-n\htop(f)}\#\Per_\chi(n,M)\) remain bounded.

\subsection{Homoclinic class and global asymptotics}

For a saddle periodic point \(x\), let \(W^s(x)\) and \(W^u(x)\)
be its global stable and unstable manifolds, and set
\(W^\sigma(\mathcal O)=\bigcup_{x\in\mathcal O}W^\sigma(x)\)
for \(\sigma\in\{s,u\}\).
Two saddle periodic orbits \(\mathcal O\) and \(\mathcal O'\)
are homoclinically related, written \(\mathcal O\sim_h\mathcal O'\),
if \(W^u(\mathcal O)\) intersects \(W^s(\mathcal O')\)
transversely and \(W^s(\mathcal O)\) intersects
\(W^u(\mathcal O')\) transversely.
The closed homoclinic class and homoclinic period of
\(\mathcal O\) are
\[
 \HH(\mathcal O)
 =\overline{\bigcup_{\mathcal O'\sim_h\mathcal O}\mathcal O'},
 \qquad
 \ell(\mathcal O)=\gcd\{\#\mathcal O':\mathcal O'\sim_h\mathcal O\}.
\]
The set \(\Per_\chi(n,\HH(\mathcal O))\) includes every
\(\chi\)-hyperbolic point of least period \(n\) in the closed
class, whether or not its orbit is homoclinically related to
\(\mathcal O\).

Buzzi, Crovisier and Sarig \cite[Theorems~1 and~2]{BCSMME}
proved a spectral decomposition at positive entropy: every
positive-entropy homoclinic class has a unique maximal measure,
whose support is the whole class, and only finitely many classes
carry global maximal measures.
They subsequently established strong positive recurrence (SPR) for
these classes \cite{BCSSPR}, linking concentration on Pesin
blocks with strongly positive recurrent countable-state
symbolic models.

In \cite{LT} we established the precise
asymptotic on each positive-entropy homoclinic class, with any
unstable cutoff below its entropy and one suitably chosen
positive stable cutoff. We also obtained a global version for transitive
diffeomorphisms.
The proof combines the recurrence theory with injective coding
and periodic-orbit estimates on the symbolic model.
Building on the estimates in \cite{LT}, we extend the range
of exponent cutoffs and remove transitivity from the global formula.
Our first result allows every common cutoff below the entropy
of the homoclinic class.

\begin{theorem}[Periodic orbit asymptotics]\label{thm:prime}
Let \(f\) be a \(C^\infty\) diffeomorphism of a closed
surface \(M\), and let \(\HH\) be a closed homoclinic class
with entropy \(h=\htop(f|\HH)>0\)
and homoclinic period \(\ell\). Then
\begin{equation}\label{eq:prime}
 \forall\,0<\chi<h,\qquad
 \lim_{\substack{n\to\infty\\\ell\mid n}}
 \frac{\#\Per_\chi(n,\HH)}{e^{nh}}=\ell.
\end{equation}
\end{theorem}

To obtain a global formula, we must also control periodic
points outside the classes supporting global maximal measures.
Proposition~\ref{prop:global-remainder} supplies the required
exponentially smaller remainder.

\Needspace{12\baselineskip}
\begin{corollary}[Global periodic counts]\label{cor:global}
Let \(f\) be a \(C^\infty\) diffeomorphism of a closed
surface \(M\) with entropy \(h=\htop(f)>0\).
Let \(\mu_1,\ldots,\mu_N\) be its
ergodic maximal measures, and let \(\ell_i\) be the homoclinic
period of \(\HH_i=\operatorname{supp}\mu_i\). Then
\begin{equation}\label{eq:global-asymptotic}
 \forall\,0<\chi<h,\qquad
 \lim_{n\to\infty}
 \left(\frac{\#\Per_\chi(n,M)}{e^{nh}}
       -\sum_{\ell_i\mid n}\ell_i\right)=0.
\end{equation}
If \(f\) is
topologically transitive, then \(N=1\).
If \(f\) is topologically mixing, then
\begin{equation}\label{eq:global-mixing}
 \forall\,0<\chi<h,\qquad
 \lim_{n\to\infty}\frac{\#\Per_\chi(n,M)}{e^{nh}}=1.
\end{equation}
\end{corollary}

\begin{remark}
The methods developed here also apply to \(C^\infty\) flows on
three-dimensional manifolds. In a forthcoming paper, we will
establish analogous results in this setting.
\end{remark}

\subsection{Lyapunov concentration}

The extension of the cutoff range rests on exponential
concentration of Lyapunov exponents.
Liao, Sun and Tian \cite{LST} identified the entropy of a
hyperbolic ergodic measure with the exponential growth rate of
periodic points approximating both the measure and its Lyapunov
exponents. For a fixed hyperbolic ergodic measure, Liao, Sun and
Yang \cite{LSY} obtained the corresponding orbit measure
approximation formula for all invariant measures carried by a
suitable set of full measure.
For expansive systems with specification, Kifer \cite{Kifer}
established large deviation estimates for periodic orbit distributions.

Our estimate concerns all \(\chi\)-hyperbolic saddle periodic points
on a compact invariant set with a unique maximal measure:
those whose exponents stay a fixed distance from the exponents
of that measure have strictly smaller exponential growth.
A unique maximal measure is ergodic and, at positive entropy,
is of saddle type by Ruelle's inequality \cite{Ruelle}.

\begin{theorem}[Exponential concentration]\label{thm:concentration}
Let \(f\) be a \(C^\infty\) diffeomorphism of a closed
surface \(M\), and let \(K\subset M\) be compact with
\(f(K)=K\) and \(h=\htop(f|K)>0\).
Suppose that \(f|K\) has a unique
maximal measure \(\mu\).
For every \(\chi,\Delta>0\), there exists \(\eta>0\) such that,
for all sufficiently large \(n\),
\begin{equation}\label{eq:concentration}
 \#\{x\in\Per_\chi(n,K):
       \max_{\sigma\in\{u,s\}}
       |\lambda^\sigma(x)-\lambda^\sigma(\mu)|\geqs\Delta\}
 \leqs e^{(h-\eta)n}.
\end{equation}
\end{theorem}

On a positive-entropy homoclinic class with maximal measure \(\mu\),
the unstable and stable cutoffs may be chosen independently in
\((0,\lambda^u(\mu))\) and \((0,-\lambda^s(\mu))\), respectively.
Theorem~\ref{thm:concentration} shows that changing either cutoff
within the corresponding range changes the count in \cite{LT} only by an
exponentially smaller term; see Corollary~\ref{cor:cutoffs}.
Ruelle's inequality then gives the range \(0<\chi<h\) in
Theorem~\ref{thm:prime}.

\subsection{Dynamical zeta functions}
For uniformly hyperbolic systems, Markov partitions provide
finite-state symbolic codings, yielding rational dynamical
zeta functions \cite{BowenLanford,Manning}. In the setting
of nonuniform hyperbolicity, however, such analytic properties
may fail, as illustrated by Kaloshin's examples of
superexponential growth of periodic points \cite{Kaloshin}.
For \(C^\infty\) surface diffeomorphisms, Burguet \cite{Burguet}
established the exponential growth rate of saddle periodic
points whose Lyapunov exponents are bounded away from zero,
thereby determining the radius of convergence of the
corresponding zeta function. In \cite[Question~7.9]{BuzziICM},
Buzzi asked whether the hyperbolic periodic-point zeta function
admits a meromorphic continuation beyond its natural disk of
convergence. An exponential refinement of our counting formula
answers this question and determines the leading poles.

Here periodic points of every hyperbolic type are included.
For \(\chi>0\), put
\[
 \Fix_\chi(n,M)=\{x\in M:f^nx=x,\ |\lambda_j(x)|>\chi
                       \text{ for }j=1,2\},
\]
where \(\lambda_1(x)\leqs\lambda_2(x)\) are the Lyapunov
exponents of the periodic orbit of \(x\).
Thus \(\Fix_\chi(n,M)\) includes sinks and sources
with least periods dividing \(n\).
Let \(\mathcal O_\chi\) denote the set of corresponding primitive
periodic orbits. The Artin--Mazur zeta function is defined by
\begin{equation}\label{eq:zeta-definition}
 \zeta_f^\chi(z)
 =\exp\left(\sum_{n\geqs1}\frac{\#\Fix_\chi(n,M)}n z^n\right)
 =\prod_{\mathcal O\in\mathcal O_\chi}
             (1-z^{\#\mathcal O})^{-1}.
\end{equation}

\begin{theorem}[Exponential counts and meromorphic continuation]
\label{thm:zeta}
Let \(f\) be a \(C^\infty\) diffeomorphism of a closed
surface with \(h=\htop(f)>0\). Let \(\mu_1,\ldots,\mu_N\)
be its ergodic maximal measures, with homoclinic periods
\(\ell_1,\ldots,\ell_N\).
For every \(0<\chi<h\), there is \(\eta>0\) such that
\begin{equation}\label{eq:hyperbolic-exponential-remainder}
 \#\Fix_\chi(n,M)
 =e^{hn}\sum_{i=1}^N\ell_i\one_{\{\ell_i\mid n\}}
       +O(e^{(h-\eta)n}).
\end{equation}
Moreover,
\begin{equation}\label{eq:hyperbolic-zeta-factorization}
 \zeta_f^\chi(z)
 =\frac{\exp H_\chi(z)}
       {\displaystyle\prod_{i=1}^N(1-e^{h\ell_i}z^{\ell_i})},
\end{equation}
where \(H_\chi\) is holomorphic on \(|z|<e^{-h+\eta}\).
The only poles in this disk are the points \(e^{-h}\omega\)
with \(\omega^{\ell_i}=1\) for some \(i\), and their orders are
\(\#\{i:\omega^{\ell_i}=1\}\).
In particular, \(e^{-h}\) is a pole of order \(N\).
\end{theorem}

Proposition~\ref{prop:class-zeta-remainder} gives an exponential
remainder for saddle counts on every positive-entropy closed
homoclinic class.
 Sinks and sources with a fixed positive exponent cutoff
have strictly smaller growth and contribute only a
holomorphic, nonvanishing factor near the leading poles.

\subsection{Proof mechanism}

Consider a homoclinic class with a unique maximal measure \(\mu\)
and entropy \(h>0\), and fix a positive hyperbolicity cutoff.
At their reference cutoffs, ~\cite[\S\S4.3--5]{LT}
handles periodic points whose orbit measures are close to \(\mu\)
and whose Lyapunov exponents are close to its exponents.
Burguet's local periodic growth estimate~\cite{Burguet}, together
with the uniqueness of \(\mu\), gives an exponential gap for periodic
points whose orbit measures stay away from \(\mu\).
It remains to obtain such a gap when orbit measures approach
\(\mu\) but the exponents retain a fixed discontinuity.

We adapt the neutral-block mechanism of Buzzi, Crovisier and
Sarig~\cite{BCSC} to periodic counting.
Lift each saddle orbit by its unstable direction.
When orbit measures approach \(\mu\), every limiting lift is a
mixture of its stable and unstable lifts. The unstable exponent
is the integral of a continuous function on the projective bundle,
so a fixed exponent defect requires positive stable mass.
The stable direction contracts for finite steps on a ``good'' block of
measure close to one. By continuity on the projective bundle,
nearby unstable directions also contract over those blocks for finite steps,
despite expanding over their full periods.

As for the direction normal to the unstable line, a hyperbolic-time type choice of origin bounds the forward normal derivatives. For parametrizations already controlled at the previous good
block, this normal bound, together with the directional contraction, controls the full endpoint derivative.
Yomdin--Gromov reparametrization and smoothness then make the entropy of the neutral part on these contracting blocks arbitrarily
small per unit time. Burguet's local estimate bounds the remaining
periodic multiplicity within each small dynamical ball.

The retained part pays the usual entropy cost,
while the neutral part can be recovered at a much smaller cost.
For finite invariant families \(\mathcal F_n\) of such points
of least period \(n\), with counting measures converging to
\(\mu\) and unstable lifts of limiting stable mass \(\omega\),
this gives
\[
 \limsup_{n\to\infty}\frac1n\log\#\mathcal F_n
 \leqs(1-\omega)h.
\]
Thus maximal growth forces \(\omega=0\).
Applied to families with a fixed exponent discrepancy, this yields
Theorem~\ref{thm:prime}. As for global counting, we use the technique of weighted Pesin blocks developed in \cite{LTTY,LT}.

\subsection{Organization of the paper}

Section~\ref{sec:prelim} develops the projective
counting estimates and Lyapunov concentration.
Section~\ref{sec:proofs} proves the main counting
estimate by erasing selected orbit names.
Section~\ref{sec:geometry} establishes the geometric
completion bound using Yomdin--Gromov reparametrization.
Section~\ref{sec:global} derives the global periodic orbit
asymptotics through homoclinic localization.
Finally, Section~\ref{sec:zeta} proves the exponential
counting remainders and the meromorphic continuation
of the hyperbolic periodic orbit zeta function.

\subsection*{Acknowledgment}
YT thanks Professor Sarig for helpful discussions in Beijing in July 2025, particularly for his explanations of the work \cite{BCSC}, which inspired some of the ideas developed in this paper. GL was partially supported by the National Key R\&D Program of China (2022YFA1005802), NSFC(12671228).
YT was partially supported by the National Key R\&D Program of China (2022YFA1005801).

\section{Reduction to a counting estimate}\label{sec:prelim}

We first use local periodic growth to identify limits of
families of maximal growth. We then formulate the counting
estimate and deduce concentration and the class asymptotics.

We use natural logarithms, with \(\log0=-\infty\) and
\(0\log0=0\) in entropy sums.
We write \(\one_A\) for the indicator of a set \(A\).
All time indices are integers. Norms are induced by the Riemannian
metric or the Euclidean norm; derivatives carry the corresponding
operator norms, and \(\|\cdot\|_\infty\) denotes the supremum norm.
All measures are Borel, and convergence of finite measures
is weak-\(*\). We write \(T_*\rho\) for the pushforward of
\(\rho\) under a Borel map \(T\).
On each compact metric space used below, \(d_{\mathrm w}\) denotes
a metric inducing the weak-\(*\) topology on probability measures.
For a nonempty finite subset \(F\) of a metric space, its normalized
counting measure is
\[
 \nu_F=\frac1{\#F}\sum_{x\in F}\delta_x,
\]
where \(\delta_x\) is the point mass at \(x\).
For \(x\in\Fix(f^n)\), its periodic orbit measure is
\[
 \nu_x=\frac1n\sum_{j=0}^{n-1}\delta_{f^jx}.
\]
This formula is unchanged if \(n\) is replaced by the least
period of \(x\). If \(F\subset\Fix(f^n)\) is invariant, then
\[
 \frac1{\#F}\sum_{x\in F}\nu_x=\nu_F.
\]
Periodic families may be indexed by any sequence of periods
tending to infinity.

All partitions below are finite Borel partitions.
For such a partition \(\PP\), let \(\PP(x)\) be its atom
containing \(x\), and let \(\partial\PP=\bigcup_{P\in\PP}\partial P\).
For \(m\geqs1\), the partition into names of length \(m\) is
\[
 \PP^m=\left\{\bigcap_{j=0}^{m-1}f^{-j}P_j:
              P_0,\ldots,P_{m-1}\in\PP\right\}
              \setminus\{\varnothing\}.
\]
For any Borel probability measure \(\rho\), put
\[
 H_\rho(\PP)=-\sum_{P\in\PP}\rho(P)\log\rho(P).
\]
If \(\rho\) is invariant, then
\(h_\rho(f,\PP)=\lim_{m\to\infty}m^{-1}H_\rho(\PP^m)\)
and \(h_\rho(f)=\sup_\PP h_\rho(f,\PP)\).

For a probability law \(\kappa\) on a finite set, write
\(H(\kappa)=-\sum_a\kappa(a)\log\kappa(a)\).
For finite-valued measurable maps \(A,B\) on a probability space with
measure \(\rho\), put
\[
 H_\rho(A)=H(A_*\rho),\qquad
 H_\rho(A\mid B)=H_\rho(A,B)-H_\rho(B).
\]
We identify a finite partition with its atom map.
For a map \(A:F\to Y\) on a nonempty finite set \(F\),
the entropy chain rule gives the counting identity
\begin{equation}\label{eq:name-decomposition}
 \log\#F=H_{\nu_F}(A)
 +\sum_{a\in A(F)}\nu_F(A^{-1}(a))\log\#A^{-1}(a).
\end{equation}

\subsection{Local periodic growth}

For each \(n\geqs1\) and \(\chi>0\), the set \(\Per_\chi(n,M)\) is finite.
Indeed, an accumulation point would still be a fixed point of
\(f^n\) with its two multipliers separated from the unit circle by
the fixed exponent bounds. It would be an isolated fixed point,
a contradiction.

For \(n\geqs1\) and \(\epsilon>0\), the \(n\)-step Bowen ball is
\[
 B_n(x,\epsilon)=
 \{y:d_M(f^jx,f^jy)<\epsilon\text{ for }0\leqs j<n\}.
\]
For a \(C^\infty\) surface diffeomorphism,
Burguet's theorem \cite[Theorem~3.2]{Burguet} gives, for every
\(\chi>0\),
\begin{equation}\label{eq:local}
 \lim_{\epsilon\downarrow0}\limsup_{n\to\infty}
 \frac1n\log\max\left\{1,\sup_{x\in M}
       \#\bigl(B_n(x,\epsilon)\cap\Per_\chi(n,M)\bigr)\right\}=0.
\end{equation}
Burguet's theorem counts all \(\chi\)-hyperbolic saddle fixed points
of \(f^n\), including those of smaller least period, so it also
applies to \(\Per_\chi(n,M)\).

The next lemma applies the argument of
\cite[Lemma~2.3]{Burguet} to varying invariant subsets.

\begin{lemma}[Periodic growth and entropy]\label{lem:periodic-entropy}
Fix \(\chi>0\). Let \(K\) be compact and invariant, and let
\(\F_n\subset\Per_\chi(n,K)\) be nonempty finite invariant sets.
If \(\nu_{\F_n}\to\nu\), then
\[
 \limsup_{n\to\infty}\frac1n\log\#\F_n\leqs h_\nu(f).
\]
In particular, their upper exponential growth rate is at most
\(\htop(f|K)\). If \(f|K\) has a unique maximal measure \(\mu\)
and \(n^{-1}\log\#\F_n\to\htop(f|K)\), then \(\nu_{\F_n}\to\mu\).
\end{lemma}

\begin{proof}
Given \(\varepsilon>0\), choose a sufficiently fine finite
partition \(\PP\) with \(\nu(\partial\PP)=0\). By \eqref{eq:local},
each atom of \(\PP^n\) contains at most
\(e^{(\varepsilon+o(1))n}\) points of \(\F_n\).
Applying \eqref{eq:name-decomposition} to the name partition,
then using invariance of \(\nu_{\F_n}\) and subadditivity of
partition entropy, gives, for each fixed \(q\),
\[
 \frac1n\log\#\F_n
 \leqs\frac1n H_{\nu_{\F_n}}(\PP^n)+\varepsilon+o(1)
 \leqs\frac1q H_{\nu_{\F_n}}(\PP^q)+\varepsilon+o(1).
\]
Let \(n\to\infty\), then \(q\to\infty\), and finally
\(\varepsilon\downarrow0\).
Compactness and the variational principle give the remaining
claims, see also \cite[Chapter~8]{Walters}.
\end{proof}

\subsection{Concentration of orbit measures}

The following consequence of ergodicity turns convergence of
counting measures into concentration of individual orbit measures.

\begin{lemma}[Concentration of orbit measures]\label{lem:orbit-measures}
Let \(T:X\to X\) be continuous on a compact metric space,
and let \(\F_n\subset\Fix(T^n)\) be nonempty finite invariant
sets. For \(x\in\F_n\), write
\(\nu_x=n^{-1}\sum_{j=0}^{n-1}\delta_{T^jx}\).
Suppose that their normalized counting measures converge to
an ergodic invariant probability measure \(\mu\).
For every \(\varepsilon>0\),
\[
 \frac{\#\{x\in\F_n:d_{\mathrm w}(\nu_x,\mu)
                         \geqs\varepsilon\}}{\#\F_n}
 \longrightarrow0.
\]
\end{lemma}

\begin{proof}
The space \(\mathcal M_T(X)\) of invariant probability measures
is compact. On this space, consider
\[
 \Theta_n=\frac1{\#\F_n}\sum_{x\in\F_n}\delta_{\nu_x}.
\]
Every subsequential limit \(\Theta\) has barycenter \(\mu\),
since \(\int\rho\,d\Theta_n(\rho)=\nu_{\F_n}\).
Ergodicity makes \(\mu\) an extreme point of
\(\mathcal M_T(X)\), so \(\Theta=\delta_\mu\).
Thus \(\Theta_n\to\delta_\mu\), and the Portmanteau theorem
applied to the closed set
\(\{\rho:d_{\mathrm w}(\rho,\mu)\geqs\varepsilon\}\)
gives the assertion.
\end{proof}

\subsection{Projective limits and the counting estimate}

For a positive-entropy ergodic measure \(\mu\), Ruelle's inequality
\cite{Ruelle}, applied to \(f\) and \(f^{-1}\), gives
\begin{equation}\label{eq:ruelle}
 \lambda^s(\mu)\leqs-h_\mu(f)<0<h_\mu(f)\leqs\lambda^u(\mu).
\end{equation}
Thus every such measure is hyperbolic of saddle type.

We now introduce the projective tangent bundle:
\[
 \wh M=\mathbb P(TM)
 =\{(x,L):x\in M,\ L\subset T_xM,\ \dim L=1\},
 \qquad \pi(x,L)=x.
\]
The induced map is
\(\wh f(x,L)=(f(x),D_xf\,L)\). Write \(J=\log|\det Df|\).
For an ergodic hyperbolic measure \(\mu\) of saddle type, fix
an invariant Borel set \(\mathcal R_\mu\) of full \(\mu\)-measure
consisting of Oseledets regular points with exponents
\(\lambda^s(\mu)<0<\lambda^u(\mu)\), at which the global Pesin stable
and unstable manifolds \(W^s(x),W^u(x)\) are defined. Write
\(T_xM=E^u(x)\oplus E^s(x)\) for the Oseledets splitting on
\(\mathcal R_\mu\); see \cite[\S2.3]{BCSMME}.
A probability measure \(\wh\rho\) on \(\wh M\) is a lift of
\(\rho\) if \(\pi_*\wh\rho=\rho\). In particular, put
\[
 \wh\mu^\sigma=(x\mapsto(x,E^\sigma(x)))_*\mu,
 \qquad \sigma\in\{u,s\}.
\]
These measures give full measure to the respective Oseledets graphs
\(\{(x,E^\sigma(x)):x\in\mathcal R_\mu\}\).
By \cite[Lemma~3.3 and Corollary~3.4]{BCSC},
every \(\wh f\)-invariant probability measure projecting to \(\mu\)
has the form \((1-\omega)\wh\mu^u+\omega\wh\mu^s\),
with \(\omega\in[0,1]\).

For a saddle periodic point \(x\), put
\(\wh x=(x,E^u(x))\).
If \(f^nx=x\), the unstable lift of its orbit measure is
\[
 \wh\nu_x=\frac1n\sum_{j=0}^{n-1}\delta_{\wh f^j\wh x}.
\]
For a nonempty finite set \(F\) of saddle periodic points,
define the unstable lift of \(\nu_F\) by
\[
 \wh\nu_F=\frac1{\#F}\sum_{x\in F}\delta_{\wh x}.
\]
It projects to \(\nu_F\) and is \(\wh f\)-invariant when
\(F\) is invariant. Thus, for invariant families \(\F_n\) with
\(\nu_{\F_n}\to\mu\), every subsequential limit of
\(\wh\nu_{\F_n}\) has the form described above.

\begin{proposition}[Periodic growth and projective limits]
\label{prop:counting}
Let \(\chi>0\), and let
\(\F_n\subset\Per_\chi(n,M)\) be nonempty invariant finite sets.
Suppose that
\[
 \nu_{\F_n}\longrightarrow\mu,
\]
where \(\mu\) is ergodic and \(h_0=h_\mu(f)>0\). Suppose also that
\[
 \wh\nu_{\F_n}\longrightarrow\wh\nu
 =(1-\omega)\wh\mu^u+\omega\wh\mu^s,
 \qquad\omega\in[0,1].
\]
Then
\begin{equation}\label{eq:counting}
 \limsup_{n\to\infty}\frac1n\log\#\F_n
 \leqs(1-\omega)h_0.
\end{equation}
\end{proposition}

The weight \(\omega\) is the mass assigned by the limiting
lift to the stable Oseledets graph.
We prove the proposition in Section~\ref{sec:proofs}.

At maximal growth, this proposition gives concentration of
individual orbit measures and Lyapunov exponents.

\begin{corollary}[Individual orbits at maximal growth]
\label{cor:maximal-growth}
Fix \(\chi>0\), and let
\(\F_n\subset\Per_\chi(n,M)\) be nonempty finite invariant sets.
Suppose that \(\nu_{\F_n}\to\mu\), where \(\mu\) is ergodic,
\(h_\mu(f)>0\), and
\[
 \frac1n\log\#\F_n\longrightarrow h_\mu(f).
\]
For every \(\varepsilon>0\),
\[
 \frac1{\#\F_n}
 \#\left\{x\in\F_n:
 d_{\mathrm w}(\nu_x,\mu)
 +\max_{\sigma\in\{u,s\}}
       |\lambda^\sigma(x)-\lambda^\sigma(\mu)|
       \geqs\varepsilon\right\}
 \longrightarrow0.
\]
\end{corollary}

\begin{proof}
Proposition~\ref{prop:counting} forces \(\omega=0\) in every
subsequential limit of the unstable lifts. Hence
\(\wh\nu_{\F_n}\to\wh\mu^u\).
Projection gives a measurable isomorphism from the unstable
Oseledets graph to \((M,\mu,f)\); hence \(\wh\mu^u\) is ergodic.
Lemma~\ref{lem:orbit-measures}, applied to the lifted families
\(\{\wh x:x\in\F_n\}\), therefore gives concentration of
\(\wh\nu_x\) at \(\wh\mu^u\).

Projection sends \(\wh\nu_x\) to \(\nu_x\), while integration
against the continuous functions
\[
 (x,L)\longmapsto\log\|D_xf|_L\|,
 \qquad
 (x,L)\longmapsto J(x)-\log\|D_xf|_L\|
\]
gives \(\lambda^u(x)\) and \(\lambda^s(x)\), respectively.
At \(\wh\mu^u\), the corresponding values are
\(\mu\), \(\lambda^u(\mu)\), and \(\lambda^s(\mu)\).
Continuity of these maps gives the asserted concentration.
\end{proof}

\subsection{Exponential concentration}

\begin{proof}[Proof of Theorem~\ref{thm:concentration}]
Recall that \(K\) is compact and invariant, \(\mu\) is its
unique maximal measure, and \(h=\htop(f|K)>0\).
Fix \(\chi,\Delta>0\), and put
\[
 \mathcal B_n=
 \left\{x\in\Per_\chi(n,K):
 \max_{\sigma\in\{u,s\}}
 |\lambda^\sigma(x)-\lambda^\sigma(\mu)|\geqs\Delta\right\}.
\]
These sets are invariant, and Lemma~\ref{lem:periodic-entropy}
bounds their upper exponential growth rate by \(h\).
If this rate were \(h\), there would be a subsequence with
\(\mathcal B_n\ne\varnothing\) and
\(n^{-1}\log\#\mathcal B_n\to h\).
Uniqueness of the maximal measure would give
\(\nu_{\mathcal B_n}\to\mu\).
Corollary~\ref{cor:maximal-growth}, with \(\varepsilon=\Delta\),
would then give a vanishing proportion of points with exponent
deviation at least \(\Delta\), whereas every point of
\(\mathcal B_n\) has this property, a contradiction.
Thus their upper growth rate is strictly smaller than \(h\),
which gives \eqref{eq:concentration}.
\end{proof}

\subsection{Independent exponent cutoffs}\label{subsec:cutoffs}

For an invariant set \(K\subset M\) and \(\alpha,\beta>0\), put
\[
 \Per_{\alpha,\beta}(n,K)
 =\{x\in K:\per(x)=n,\ \lambda^u(x)>\alpha,\ \lambda^s(x)<-\beta\},
\]
where only saddle points are included. Thus
\(\Per_\chi(n,K)=\Per_{\chi,\chi}(n,K)\).
By \eqref{eq:ruelle}, the following consequence of concentration
includes Theorem~\ref{thm:prime}.

\begin{corollary}[Independent exponent cutoffs]\label{cor:cutoffs}
Let \(f\) be a \(C^\infty\) diffeomorphism of a closed
surface \(M\), and let \(\HH\) be a closed homoclinic class
with entropy \(h=\htop(f|\HH)>0\), homoclinic period
\(\ell\), and maximal measure \(\mu\). Then
\begin{equation}\label{eq:independent-cutoffs}
 \forall\,(\alpha,\beta)\in
 (0,\lambda^u(\mu))\times(0,-\lambda^s(\mu)),\qquad
 \lim_{\substack{n\to\infty\\\ell\mid n}}
 \frac{\#\Per_{\alpha,\beta}(n,\HH)}{e^{nh}}=\ell.
\end{equation}
\end{corollary}

\begin{proof}[Proof of Corollary~\ref{cor:cutoffs} and
Theorem~\ref{thm:prime}]
The class \(\HH\) has a unique maximal measure \(\mu\) by
\cite[Theorem~2]{BCSMME}.
Fix \(\alpha_*=h/2\).
By \cite[Theorem~1.1]{LT},
there is a cutoff \(\beta_*\in(0,-\lambda^s(\mu))\) such that
\begin{equation}\label{eq:LT}
 \lim_{\substack{n\to\infty\\ \ell\mid n}}
 \frac{\#\Per_{\alpha_*,\beta_*}(n,\HH)}{e^{nh}}=\ell.
\end{equation}
The signed stable cutoff in \cite[\S4.3, equation~(4)]{LT}
lies in \((\lambda^s(\mu),0)\), which gives the stated range
for \(\beta_*\).

Fix \(0<\alpha<\lambda^u(\mu)\) and \(0<\beta<-\lambda^s(\mu)\).
Every point in the symmetric difference
\(\Per_{\alpha,\beta}(n,\HH)\triangle
  \Per_{\alpha_*,\beta_*}(n,\HH)\)
is \(\min\{\alpha,\beta,\alpha_*,\beta_*\}\)-hyperbolic and satisfies
\[
 \begin{aligned}
 &\lambda^u(x)\leqs\max\{\alpha,\alpha_*\}<\lambda^u(\mu)
 \quad\text{or}\\
 &\lambda^s(x)\geqs-\max\{\beta,\beta_*\}>\lambda^s(\mu).
 \end{aligned}
\]
Since the cutoffs are fixed, Theorem~\ref{thm:concentration}
gives, for some \(\eta>0\) and all large \(n\),
\begin{equation}\label{eq:cutoff-comparison}
 \#\bigl(\Per_{\alpha,\beta}(n,\HH)
       \triangle\Per_{\alpha_*,\beta_*}(n,\HH)\bigr)
 \leqs e^{(h-\eta)n}.
\end{equation}
Together with \eqref{eq:LT}, this proves
Corollary~\ref{cor:cutoffs}.
Finally, set \(\alpha=\beta=\chi<h\) and use \eqref{eq:ruelle}
to obtain Theorem~\ref{thm:prime}.
\end{proof}

\section{Erasing orbit names and the counting estimate}\label{sec:proofs}

We prove Proposition~\ref{prop:counting} by erasing intervals
for which the decrease in the retained-name entropy bound
exceeds the cost of completion.
The separation of these two contributions follows the
neutral-block strategy of \cite[\S\S2.2--2.4]{BCSC};
see also \cite[\S\S2, 5--6]{BuzziLectures}.
Here this strategy is applied to finite periodic families.
We first establish the normal derivative bound used in the
completion estimate.

\subsection{A cyclic origin with bounded normal derivative}
\label{subsec:normal}

Along a saddle periodic orbit, the derivative induced on the
quotient by the unstable line is one-dimensional, and the
logarithm of its norm has a negative sum over a period.
Choosing a cyclic origin at a maximum partial sum therefore
gives uniform normal control throughout the period.
For a line \(L\), let \(\pi_L\) denote orthogonal projection
onto \(L\).

\begin{lemma}[Bounded normal derivatives]\label{lem:normal}
For every saddle point \(x\) fixed by \(f^n\), there is
\(\tau\in\{0,\ldots,n-1\}\) such that, with \(y=f^\tau x\),
\[
 \|\pi_{E^u(f^ty)^\perp}D_yf^t\|\leqs1
 \qquad(0\leqs t\leqs n).
\]
\end{lemma}

\begin{proof}
Along the orbit of \(x\), put
\[
 a(z)=\log|\det D_zf|-\log\|D_zf|_{E^u(z)}\|.
\]
For \(A=D_zf^t\), choose orthonormal bases with first vectors
in \(E^u(z)\) and \(E^u(f^tz)\). The matrix of \(A\) is
upper triangular, with second row \((0,d)\) and
\(|d|=|\det A|/\|A|_{E^u(z)}\|\). Thus
\begin{equation}\label{eq:normal-area}
 \|\pi_{E^u(f^tz)^\perp}D_zf^t\|
 =\frac{|\det D_zf^t|}{\|D_zf^t|_{E^u(z)}\|}
 =\exp\left(\sum_{j=0}^{t-1}a(f^jz)\right).
\end{equation}
Set \(S_t=\sum_{j=0}^{t-1}a(f^jx)\). Since
\(S_n=n\lambda^s(x)<0\), an index \(\tau\) maximizing
\(S_t\) over \(0\leqs t<n\) satisfies
\[
 \sum_{j=0}^{t-1}a(f^{\tau+j}x)
 =\begin{cases}
 S_{\tau+t}-S_\tau,&\tau+t\leqs n,\\
 S_n-S_\tau+S_{\tau+t-n},&\tau+t>n,
 \end{cases}
 \leqs0
 \qquad(0\leqs t\leqs n).
\]
Equation~\eqref{eq:normal-area} now gives the assertion
with \(y=f^\tau x\).
\end{proof}

Fix a smooth atlas with \(C_0\) charts
\(\varphi_i:U_i\to(-3,3)^2\), such that the sets
\(U_i'=\varphi_i^{-1}(B(0,1))\) cover \(M\), and put
\begin{equation}\label{eq:coordinate-maps}
 \Phi_i=\varphi_i^{-1}|_{(-2,2)^2}.
\end{equation}
Fix \(C_2\geqs1\) such that every
\(\varphi_i^{-1}|_{[-2,2]^2}\) is \(C_2\)-bi-Lipschitz.
For \(y=f^\tau x\) supplied by Lemma~\ref{lem:normal},
choose \(i\) with \(y\in U_i'\) and put \(u=\varphi_i(y)\).
The chain rule gives
\begin{equation}\label{eq:normal}
 \|\pi_{E^u(f^ty)^\perp}D_u(f^t\circ\Phi_i)\|\leqs C_2
 \qquad(0\leqs t\leqs n).
\end{equation}
The constant \(C_2\) depends only on the metric and the atlas.

\subsection{Entropy of retained names}

For the rest of this section and Section~\ref{sec:geometry},
assume the hypotheses of Proposition~\ref{prop:counting}.
Write
\[
 h_0=h_\mu(f)>0,\qquad
 \wh\nu=(1-\omega)\wh\mu^u+\omega\wh\mu^s.
\]
For each \(x\in\F_n\), choose a starting time \(\tau(x)\)
from Lemma~\ref{lem:normal}, then choose \(i=i(x)\) with
\(f^{\tau(x)}x\in U_i'\) and write \(\Phi_{0,x}=\Phi_i\).
Put
\begin{equation}\label{eq:transported-lifts}
 \begin{aligned}
 u_x&=\Phi_{0,x}^{-1}(f^{\tau(x)}x),\qquad
 L_t(x)=E^u(f^{\tau(x)+t}x),\\
 \wh z_t(x)&=(f^{\tau(x)+t}x,L_t(x)).
 \end{aligned}
\end{equation}
These lifted points satisfy
\(\wh z_t(x)=\wh f^t\wh z_0(x)\).
There are at most \(C_0n\) pairs \((\tau(x),\Phi_{0,x})\).
Complete cyclic sums and the invariance of \(\F_n\) give
\begin{equation}\label{eq:transported-limit}
 \frac1{n\#\F_n}\sum_{x\in\F_n}\sum_{t=0}^{n-1}
             \delta_{\wh z_t(x)}
 =\wh\nu_{\F_n}\longrightarrow\wh\nu.
\end{equation}

For a finite partition \(\PP\) and arbitrary sets of erased
times \(I_n(x)\subset\{0,\ldots,n-1\}\), \(x\in\F_n\), put
\(R_n(x)=\{0,\ldots,n-1\}\setminus I_n(x)\).
Adjoin one symbol \(\ast\) to \(\PP\), and put
\(\mathcal A=\PP\sqcup\{\ast\}\).
Define the word map \(W_n:\F_n\to\mathcal A^n\) by
\[
 W_n(x)_j=
 \begin{cases}
 \PP(f^{\tau(x)+j}x),&j\in R_n(x),\\
 \ast,&j\in I_n(x),
 \end{cases}
 \quad 0\leqs j<n.
\]
The partial name map on \(\F_n\) is
\begin{equation}\label{eq:partial-name-map}
 C(x)=\bigl(\tau(x),\Phi_{0,x},W_n(x)\bigr).
\end{equation}
Since \(R_n(x)=\{j:W_n(x)_j\ne\ast\}\), this name determines
the retained times and all their symbols.
Equation~\eqref{eq:name-decomposition} separates the entropy
of this name from the average logarithmic number of its
completions.

The following estimate adapts Burguet's partial-name argument
\cite[Lemma~6]{BurguetSRB}, \cite[Lemma~4]{BurguetContinuity}
to cyclic words. When the number of erased intervals is small,
most short windows in \(W_n\) are wholly erased or wholly retained.
Cyclic averaging then relates their entropy to the occupation
measure on the retained times.

\begin{lemma}[Retained-name entropy]\label{lem:retained}
Fix a finite partition \(\PP\) with \(\mu(\partial\PP)=0\).
For every \(\varepsilon>0\), there are \(\beta>0\) and \(n_0\geqs1\)
such that, for every \(n\geqs n_0\), if each \(I_n(x)\)
is a union of at most \(\beta n\) intervals of consecutive
integers, then
\begin{equation}\label{eq:retained}
 H_{\nu_{\F_n}}(C)
 \leqs h_0\int|R_n|\,d\nu_{\F_n}+\varepsilon n.
\end{equation}
The estimate is uniform over all such choices of \(I_n\).
\end{lemma}

\begin{proof}
Suppose the conclusion fails for some \(\varepsilon>0\).
Along a sequence of periods \(n\to\infty\), there are then
integers \(b_n=o(n)\) and choices of erased sets, each a union
of at most \(b_n\) intervals, for which \eqref{eq:retained}
fails. Define
\[
 \eta_n=\frac1{n\#\F_n}
       \sum_{x\in\F_n}\sum_{j\in R_n(x)}
       \delta_{f^{\tau(x)+j}x},
 \qquad
 t_n=\eta_n(M)=\frac1n\int|R_n|\,d\nu_{\F_n}.
\]
Averaging complete cyclic sums and using invariance of
\(\F_n\) gives \(\eta_n\leqs\nu_{\F_n}\).
Cancellation at the endpoints of the at most \(b_n+1\)
retained intervals gives
\[
 \|f_*\eta_n-\eta_n\|_{\mathrm{TV}}
 \leqs\frac{2(b_n+1)}n.
\]
Passing to a subsequence, assume \(\eta_n\to\eta\) and
\(t_n\to t\in[0,1]\). The preceding estimates give
\(f_*\eta=\eta\) and \(\eta\leqs\mu\).
Ergodicity of \(\mu\) therefore gives \(\eta=t\mu\).

Write \(W_j(x)=W_n(x)_j\). For \(1\leqs m\leqs n\), define
the probability measure on \(\mathcal A^m\)
\[
 \kappa_{n,m}
 =\frac1{n\#\F_n}
   \sum_{x\in\F_n}\sum_{s=0}^{n-1}
   \delta_{\bigl(W_{(s+j)\bmod n}(x)\bigr)_{0\leqs j<m}}.
\]
For \(0\leqs s<n\), put
\(J_s=\{(s+j)\bmod n:0\leqs j<m\}\).
The entropy chain rule, applied in increasing order of
the indices in \(J_s\), gives
\[
 \begin{aligned}
 H_{\nu_{\F_n}}
   \bigl((W_{(s+j)\bmod n})_{0\leqs j<m}\bigr)
 &=\sum_{i\in J_s}H_{\nu_{\F_n}}
   \bigl(W_i\mid(W_j)_{j\in J_s,\ j<i}\bigr)\\
 &\geqs\sum_{i\in J_s}H_{\nu_{\F_n}}
   (W_i\mid W_0,\ldots,W_{i-1}).
 \end{aligned}
\]
Each index belongs to exactly \(m\) of the sets \(J_s\), and
\(\kappa_{n,m}\) averages the distributions of these words
over \(s\). Summing over \(s\), using concavity of entropy
and the \(C_0n\) bound on the number of pairs
\((\tau(x),\Phi_{0,x})\), we obtain
\[
 \begin{aligned}
 H_{\nu_{\F_n}}(C)
 &\leqs\log(C_0n)+H_{\nu_{\F_n}}(W_n)\\
 &\leqs\log(C_0n)+\frac1m\sum_{s=0}^{n-1}
   H_{\nu_{\F_n}}
   \bigl((W_{(s+j)\bmod n})_{0\leqs j<m}\bigr)\\
 &\leqs\log(C_0n)+\frac nm H(\kappa_{n,m}).
 \end{aligned}
\]

Now fix \(m\). For each \(x\), any \(J_s\) meeting both
\(R_n(x)\) and \(I_n(x)\) contains an endpoint of a retained
interval. There are at most \(2(b_n+1)\) such endpoints,
each belonging to \(m\) of the sets \(J_s\).
Thus the proportion of these pairs \((x,s)\) is at most
\(2m(b_n+1)/n\to0\).

For \(w=(A_0,\ldots,A_{m-1})\) with \(A_j\in\PP\), put
\(B_w=\bigcap_{j=0}^{m-1}f^{-j}A_j\).
By definition,
\[
 \eta_n(B_w)
 =\frac1{n\#\F_n}\sum_{x\in\F_n}
   \#\{s\in R_n(x):f^{\tau(x)+s}x\in B_w\}.
\]
Except for the mixed windows counted above, the word at \(s\)
is \(\ast^m\) when \(s\in I_n(x)\), and is the ordinary
length-\(m\) \(\PP\)-name of \(f^{\tau(x)+s}x\) when
\(s\in R_n(x)\).
This remains valid when \(s+m>n\), since \(f^nx=x\).
The mixed windows have vanishing frequency, while
\(\eta_n\to t\mu\) and \(\mu(\partial B_w)=0\). Hence
\[
 \begin{aligned}
 \kappa_{n,m}(\{\ast^m\})&=1-t_n+o(1)\longrightarrow1-t,\\
 \kappa_{n,m}(\{w\})
 &=\eta_n(B_w)+o(1)\longrightarrow t\mu(B_w).
 \end{aligned}
\]
The total mass of words containing both \(\ast\) and symbols
from \(\PP\) tends to zero. The nonempty sets \(B_w\) are the
atoms of \(\PP^m\), so continuity of entropy on the fixed
finite alphabet \(\mathcal A^m\) gives
\[
 \begin{aligned}
 \lim_{n\to\infty}H(\kappa_{n,m})
 &=-(1-t)\log(1-t)-t\log t+tH_\mu(\PP^m)\\
 &\leqs\log2+tH_\mu(\PP^m).
 \end{aligned}
\]
Consequently, for every fixed \(m\),
\[
 \limsup_{n\to\infty}\frac1nH_{\nu_{\F_n}}(C)
 \leqs\frac{\log2+tH_\mu(\PP^m)}m.
\]
Letting \(m\to\infty\) and using
\(m^{-1}H_\mu(\PP^m)\to h_\mu(f,\PP)\leqs h_0\),
we obtain an upper bound of \(th_0\).
On the other hand, the failure of \eqref{eq:retained} and
\(t_n\to t\) imply
\[
 \liminf_{n\to\infty}\frac1nH_{\nu_{\F_n}}(C)
 \geqs th_0+\varepsilon,
\]
a contradiction.
\end{proof}
\subsection{Good blocks and the completion bound}\label{subsec:completion}

For \(p\geqs1\), define the set of good block starts by
\begin{equation}\label{eq:good-projective-set}
 \Gamma_p=\{(x,L)\in\wh M:\|D_xf^p|_L\|<1\}.
\end{equation}
This set is open. Since \(\lambda^s(\mu)<0\),
\(\|D_xf^p|_{E^s(x)}\|\to0\) for \(\mu\)-almost every \(x\).
Bounded convergence therefore gives
\begin{equation}\label{eq:stable-block-mass}
 \wh\mu^s(\Gamma_p)\longrightarrow1\qquad(p\to\infty).
\end{equation}

We erase complete blocks on a common grid.
Fix \(p\geqs1\), take \(n>p\), and choose an offset
\(k\in\{0,\ldots,p-1\}\). Put
\[
 m_k=\left\lfloor\frac{n-1-k}{p}\right\rfloor,\qquad
 t_j=k+jp\quad(0\leqs j\leqs m_k).
\]
For each \(x\in\F_n\), choose
\(\mathcal I(x)\subset\{1,\ldots,m_k\}\) and set
\begin{equation}\label{eq:erased-times}
 I_n(x)=\bigcup_{j\in\mathcal I(x)}
             \{t_{j-1}+1,\ldots,t_j\}.
\end{equation}
The initial times \(0,\ldots,k\) and the final times
\(t_{m_k}+1,\ldots,n-1\) are retained.
For the set \(\Gamma_p\) in \eqref{eq:good-projective-set}, define
\[
 N_n^{\mathrm{out}}(x)
 =\#\{j\in\mathcal I(x):\wh z_{t_{j-1}}(x)\notin \Gamma_p\}.
\]
We charge an additional completion cost for each such block.

\begin{proposition}[Completion of a partial name]\label{prop:completion}
With the notation above, there is a constant \(C_4>0\),
depending only on \(f\), the metric, and the fixed atlas,
with the following property.
For every \(\varepsilon>0\), there is
an integer \(p_1\geqs1\) such that, for each fixed \(p\geqs p_1\),
one can choose a finite partition \(\PP\) with
\(\mu(\partial\PP)=0\) for which
\begin{equation}\label{eq:completion-entropy}
 \log\#\F_n-H_{\nu_{\F_n}}(C)
 \leqs\varepsilon n+
       C_4p\int N_n^{\mathrm{out}}\,d\nu_{\F_n}
\end{equation}
for all sufficiently large \(n\), uniformly over all offsets
\(k\) and all choices of \(\mathcal I(x)\) in
\eqref{eq:erased-times}.
\end{proposition}

Section~\ref{sec:geometry} proves this bound by successive
reparametrizations, with local periodic growth controlling
the final multiplicity.

\subsection{Selection of erased intervals}

Grouping consecutive \(p\)-blocks keeps the number of erased
intervals small enough for Lemma~\ref{lem:retained}.
The retained-name entropy bound decreases by \(h_0p\) for
each erased \(p\)-block, while Proposition~\ref{prop:completion}
charges an additional cost \(C_4p\) when its initial lifted
point lies outside \(\Gamma_p\).
For integers \(p,q\geqs1\), define
\(\Psi_{p,q}:\wh M\to[0,h_0]\) by
\[
 \Psi_{p,q}(\wh z)
 =\max\left\{0,\,
 h_0-\frac{C_4}{q}\sum_{i=0}^{q-1}
       \one_{\Gamma_p^c}(\wh f^{ip}\wh z)\right\},
\]
where \(\Gamma_p\) is given by \eqref{eq:good-projective-set}.
Thus \(pq\,\Psi_{p,q}(\wh z)\) is the net saving in these
bounds from a group of \(q\) blocks starting at \(\wh z\), obtained by
erasing it when its saving exceeds its cost and retaining
it otherwise.

\begin{lemma}[Selection of erased intervals]\label{lem:selection}
For every \(\varepsilon>0\), there is
an integer \(p_2\geqs1\) such that, for each fixed integer
\(p\geqs p_2\), one can choose an integer \(q\geqs1\) for which
\begin{equation}\label{eq:selection-counting}
 \frac1n\log\#\F_n
 \leqs h_0-\int\Psi_{p,q}\,d\wh\nu_{\F_n}+\varepsilon
\end{equation}
holds for all sufficiently large \(n\).
\end{lemma}

\begin{proof}
Apply Proposition~\ref{prop:completion} with error
\(\varepsilon/3\), and denote the resulting threshold by
\(p_2\). Fix \(p\geqs p_2\) and the corresponding partition
\(\PP\).
Apply Lemma~\ref{lem:retained} with error \(\varepsilon/3\)
to obtain \(\beta>0\).
Choose an integer \(q\) with \((pq)^{-1}<\beta\), and put \(T=pq\).

For \(n>T\) and each offset \(s\in\{0,\ldots,T-1\}\), consider
groups of \(q\) consecutive \(p\)-blocks beginning at
\(t=s+jT\), where \(j\geqs0\) and \(t+T<n\).
For each \(x\in\F_n\), erase the times
\(t+1,\ldots,t+T\) if \(\Psi_{p,q}(\wh z_t(x))>0\), and
retain them otherwise. Retain all times outside these groups.
Since \(p\mid T\), this selection has the form
\eqref{eq:erased-times} with offset \(k=s\bmod p\).
The erased set has at most \(n/T<\beta n\) intervals.

For each \(x\in\F_n\), the resulting saving is
\[
 h_0|I_n(x)|-C_4pN_n^{\mathrm{out}}(x)
 =T\sum_{\substack{j\geqs0\\s+(j+1)T<n}}
       \Psi_{p,q}(\wh z_{s+jT}(x)).
\]
For all sufficiently large \(n\), estimates
\eqref{eq:retained} and \eqref{eq:completion-entropy} apply
simultaneously to all these choices.
Add them and average over the \(T\) offsets.
Each time \(t\in\{0,\ldots,n-T-1\}\) is a group starting
time for exactly one offset, so
\[
 \begin{aligned}
 \log\#\F_n
 &\leqs h_0n-\frac1{\#\F_n}
    \sum_{x\in\F_n}\sum_{t=0}^{n-T-1}
       \Psi_{p,q}(\wh z_t(x))+\frac{2\varepsilon}{3}n\\
 &\leqs n\left(h_0-\int\Psi_{p,q}\,d\wh\nu_{\F_n}
                   +\frac{2\varepsilon}{3}\right)+h_0T.
 \end{aligned}
\]
The last inequality follows from \eqref{eq:transported-limit}
and \(0\leqs\Psi_{p,q}\leqs h_0\).
Since \(T\) is fixed, \(h_0T\leqs\varepsilon n/3\) for all
sufficiently large \(n\). Dividing by \(n\) proves
\eqref{eq:selection-counting}.
\end{proof}

\subsection{Proof of the counting estimate}

\begin{proof}[Proof of Proposition~\ref{prop:counting}]
Let \(C_4\) be the constant in Proposition~\ref{prop:completion}
and fix \(\varepsilon>0\).
By \eqref{eq:stable-block-mass}, choose \(p\) sufficiently
large that \(\wh\mu^s(\Gamma_p^c)\leqs\varepsilon\) and
Lemma~\ref{lem:selection} applies. Fix the corresponding \(q\).
These parameters remain fixed as \(n\to\infty\).

Recall from \eqref{eq:transported-limit} that
\(\wh\nu_{\F_n}\to\wh\nu\).
Since \(\Gamma_p\) is open, \(\Psi_{p,q}\) is lower semicontinuous.
Its nonnegativity and invariance of \(\wh\mu^s\) give
\[
 \begin{aligned}
 \int\Psi_{p,q}\,d\wh\nu
 &\geqs\omega\int\Psi_{p,q}\,d\wh\mu^s\\
 &\geqs\omega\bigl(h_0-C_4\wh\mu^s(\Gamma_p^c)\bigr)
 \geqs h_0\omega-C_4\varepsilon.
 \end{aligned}
\]
Lower semicontinuity then gives
\[
 \liminf_{n\to\infty}\int\Psi_{p,q}\,d\wh\nu_{\F_n}
 \geqs h_0\omega-C_4\varepsilon.
\]
Equation~\eqref{eq:selection-counting} now gives
\[
 \begin{aligned}
 \limsup_{n\to\infty}\frac1n\log\#\F_n
 &\leqs h_0-\liminf_{n\to\infty}\int\Psi_{p,q}\,d\wh\nu_{\F_n}
          +\varepsilon\\
 &\leqs(1-\omega)h_0+(C_4+1)\varepsilon.
 \end{aligned}
\]
Since \(C_4\) is independent of \(\varepsilon\), letting
\(\varepsilon\downarrow0\) proves the proposition.
\end{proof}

\section{Proof of the completion estimate}\label{sec:geometry}

By \eqref{eq:name-decomposition}, the completion cost in
Proposition~\ref{prop:completion} is the average logarithmic
cardinality of the fibers of \(C\).
We estimate it by successive block coverings with
\(C^r\)-bounded parametrizations.

The key step is Lemma~\ref{lem:block-covering}(2):
on good erased blocks, contraction along the unstable line
and the normal bound of Lemma~\ref{lem:normal} together control
the full derivative of the endpoint parametrization at the
parameters representing periodic points.
Lemma~\ref{lem:reparametrization} then turns this pointwise
control into a covering with small exponential cost.
Retained blocks admit the same bound through their prescribed
endpoint atoms (Lemma~\ref{lem:block-covering}(1)).
Lemma~\ref{lem:finite-scale-completion} averages the additional
cost of the other erased blocks, while local periodic growth
\eqref{eq:local} controls the remaining multiplicity.

\subsection{Local reparametrization}

For an integer \(r\geqs2\), use the usual \(C^r\) derivative seminorm
\[
 |F|_r=\max_{1\leqs k\leqs r}\|D^kF\|_\infty
\]
with suprema over the specified domain.
Maps on the closed parameter square \([0,1]^2\) are understood
to extend to \(C^r\) maps on a neighbourhood.
We use the functional form of the Yomdin--Gromov
reparametrization lemma \cite{Gromov,Yomdin}, as stated in
\cite[Lemma~5]{BurguetSurface}, for both values and first derivatives.

\begin{lemma}[Local reparametrization]\label{lem:reparametrization}
For every integer \(r\geqs2\) and every \(B\geqs1\), there is a
constant \(C_1(r,B)\geqs1\) with the following property.
Let \(F:[0,1]^2\to\R^2\) be \(C^r\), with
\(\|D^rF\|_\infty\leqs A\), \(A\geqs1\).
For \(v_0\in\R^2\), set
\[
 S_1=\{u\in[0,1]^2:\|DF(u)\|\leqs B\},
 \qquad S_2=F^{-1}(\overline B(v_0,B)).
\]
For each \(j=1,2\), there are at most
\(C_1(r,B)A^{2/r}\) smooth maps
\(\theta_i:[0,1]^2\to[0,1]^2\) such that
\[
 S_j\subset\bigcup_i\theta_i([0,1]^2),\qquad
 \|D\theta_i\|_\infty\leqs1,\qquad
 |F\circ\theta_i|_r\leqs1.
\]
\end{lemma}

\begin{proof}
Fix \(j\in\{1,2\}\) and put \(m=\lceil(2A)^{1/r}\rceil\).
The factor \(2\) accounts for the change from operator to
Euclidean matrix norms.
Subdivide \([0,1]^2\) into \(m^2\) closed squares with affine
parametrizations \(\ell_c(v)=a_c+m^{-1}v\).
Set \(F_c=F\circ\ell_c\) and \(S_c=\ell_c^{-1}(S_j)\).
For the functional lemma, take
\((g,s)=(F_c-v_0,r)\) for \(j=2\), and
\((g,s)=(DF_c,r-1)\) for \(j=1\), viewing \(DF_c\) as an
\(\R^4\)-valued map.
Since \(DF_c=m^{-1}(DF)\circ\ell_c\), in both cases
\[
 \|D^sg\|_\infty\leqs\sqrt2\,m^{-r}A\leqs1,
 \qquad \|g(u)\|\leqs B\quad(u\in S_c),
\]
where the second bound uses \(m\geqs2\) in the derivative
case. In particular, the Lipschitz seminorm of \(D^{s-1}g\)
is at most one, so \cite[Lemma~5]{BurguetSurface} applies with
\(B\). It gives \(O_r(1)\) analytic maps
\(\theta:[0,1]^2\to[0,1]^2\) whose images cover \(S_c\), with
\[
 |\theta|_{s+1}\leqs1,\qquad
 \max_{1\leqs\ell\leqs s}
       \|D^\ell(g\circ\theta)\|_\infty\leqs B.
\]
In the derivative case, \(s+1=r\), so the first estimate
controls \(|\theta|_r\).

For \(j=2\), these estimates give
\(|F_c\circ\theta|_r\leqs B\).
For \(j=1\), retain only maps whose images meet \(S_c\).
Then
\[
 \|g\circ\theta\|_\infty
 \leqs B+\sqrt2\,\|D(g\circ\theta)\|_\infty
 \leqs(1+\sqrt2)B.
\]
Thus \(g\circ\theta\) and \(D\theta\) have \(C^{r-1}\) bounds
depending only on \(r,B\). Applying the product rule to
\[
 D(F_c\circ\theta)=(g\circ\theta)D\theta
\]
gives \(|F_c\circ\theta|_r\leqs C(r,B)\).

In both cases, \(|F_c\circ\theta|_r\leqs C(r,B)\).
Subdivide the domain of each \(\theta\) into
\(O_{r,B}(1)\) equal squares whose side length \(b\) satisfies
\(b\leqs\min\{1,C(r,B)^{-1}\}\).
For each affine parametrization \(\eta\) of such a square,
the chain rule gives
\[
 |F_c\circ\theta\circ\eta|_r
 \leqs b|F_c\circ\theta|_r\leqs1,
 \qquad
 \|D(\ell_c\circ\theta\circ\eta)\|_\infty
 \leqs m^{-1}b\leqs1.
\]
The maps \(\ell_c\circ\theta\circ\eta\) therefore give the
required covering by \(O_{r,B}(m^2)=O_{r,B}(A^{2/r})\) maps.
\end{proof}

\subsection{Covering one block}

\emph{Scales and partitions.}
Recall the atlas \(\varphi_i:U_i\to(-3,3)^2\) and coordinate
maps \(\Phi_i\) from \eqref{eq:coordinate-maps}.
Each \(\varphi_i^{-1}|_{[-2,2]^2}\) is \(C_2\)-bi-Lipschitz,
with \(C_2\geqs1\) fixed in Subsection~\ref{subsec:normal}.
For \(a\geqs1\), let
\[
 F_a^{ij}=\varphi_j\circ f^a\circ\varphi_i^{-1}
 \quad\text{on }\varphi_i(U_i\cap f^{-a}U_j).
\]
Derivative suprema and seminorms for these maps are taken over
\([-2,2]^2\cap(F_a^{ij})^{-1}([-2,2]^2)\), and are zero
when this set is empty.
Choose
\[
 \Lambda>C_2^2
 \max\{1,\|Df\|_\infty,\|Df^{-1}\|_\infty\}.
\]
Then, for every \(a\geqs1\),
\[
 \|Df^{\pm a}\|_\infty\leqs\Lambda^a,
 \qquad
 \max_{i,j}\|DF_a^{ij}\|_\infty
 \leqs C_2^2\|Df\|_\infty^a
 \leqs\Lambda^a.
\]
The constant \(\Lambda\) depends only on \(f\), the metric,
and the atlas, and is fixed throughout this section.

Fix integers \(r\geqs2\) and \(p\geqs1\).
Choose \(\delta_0(r,p)\in(0,1]\) so small that, whenever
\(0<\delta\leqs\delta_0(r,p)\),
\begin{equation}\label{eq:jets}
 \delta\max_{\substack{1\leqs a\leqs p\\i,j}}
       |F_a^{ij}|_r\leqs1,
\end{equation}
and every subset of \(M\) of diameter at most
\(\sqrt2 C_2\Lambda^p\delta\) is contained in some
\(\varphi_j^{-1}((-2,2)^2)\).
The derivative condition involves finitely many smooth maps
on compact sets; the geometric condition follows from the
Lebesgue number lemma for the inner chart domains.
They control the higher derivatives after rescaling and
ensure that each block image used below lies in a single chart.

For each \(0<\delta\leqs\delta_0(r,p)\), choose a finite
partition \(\PP_\delta\) with
\[
 \mu(\partial\PP_\delta)=0,\qquad
 \max_{P\in\PP_\delta}\operatorname{diam}P\leqs\delta.
\]
Such partitions are obtained from finite covers by balls of
radius less than \(\delta/2\) with \(\mu\)-null boundaries.

\emph{\(C^r\)-bounded parametrizations and normal control.}
Fix \(0<\delta\leqs\delta_0(r,p)\) and use
\(\PP=\PP_\delta\) in \eqref{eq:partial-name-map}.
For \(n>p\), fix a grid offset \(k\in\{0,\ldots,p-1\}\)
and an arbitrary
erased-block selection \(\mathcal I(x)\) as in \eqref{eq:erased-times}.
Recall that \(t_j=k+jp\) and
\(m_k=\lfloor(n-1-k)/p\rfloor\); the sampled times are
\(\mathcal T_k=\{0,t_0,\ldots,t_{m_k}\}\).
Fix \(c\in C(\F_n)\) and put
\(\F_c=C^{-1}(c)\).
On this fiber, \(\tau,\Phi_0,W_n\) are fixed, so the erased
blocks and all retained endpoint atoms are prescribed.
Write \(\Phi_0=\Phi_i\) and
\(V_t=f^t\circ\Phi_0\).
For \(x\in\F_c\), recall that
\[
 u_x=\Phi_0^{-1}(f^\tau x),\qquad
 L_t(x)=E^u(f^{\tau+t}x).
\]
The corresponding lifted point is
\(\wh z_t(x)=(f^{\tau+t}x,L_t(x))\).
The normal derivative bound \eqref{eq:normal} reads
\[
 \|\pi_{L_t(x)^\perp}D_{u_x}V_t\|\leqs C_2
 \qquad(0\leqs t\leqs n).
\]

We use smooth parameter maps
\(\sigma:[0,1]^2\to[-1,1]^2\) in the initial coordinate domain,
all satisfying
\begin{equation}\label{eq:original-coordinates}
 \|D\sigma\|_\infty\leqs\frac{\delta}{C_2\Lambda^p}.
\end{equation}
Such a map \emph{covers} \(x\in\F_c\) if
\(u_x\in\sigma([0,1]^2)\).
The parametrization in \(M\) at time \(t\) is \(V_t\circ\sigma\).
At the fixed scale \(\delta\), for a chart \(\varphi\) from the
atlas, the pair \((\sigma,\varphi)\) is
\emph{\(C^r\)-bounded at time \(t\)} if
\[
 V_t(\sigma([0,1]^2))\subset\varphi^{-1}((-2,2)^2),
 \qquad |\varphi\circ V_t\circ\sigma|_r\leqs\delta.
\]

Condition \eqref{eq:original-coordinates} is preserved by every precomposition
\(\sigma\mapsto\sigma\circ\theta\) with
\(\theta:[0,1]^2\to[0,1]^2\) and
\(\|D\theta\|_\infty\leqs1\).
For each \(x\in\F_c\) and each parameter \(v\) with
\(\sigma(v)=u_x\), the chain rule and \eqref{eq:normal} give
\begin{equation}\label{eq:projection}
 \|\pi_{L_t(x)^\perp}D_v(V_t\circ\sigma)\|
 \leqs C_2\|D_v\sigma\|\leqs\delta\Lambda^{-p}
 \qquad(0\leqs t\leqs n).
\end{equation}
The factor \(\Lambda^{-p}\) in this normal bound will
offset the inverse-derivative growth over a block of length \(p\).

\begin{lemma}[Covering one block]\label{lem:block-covering}
There is a constant \(C_3(r)\geqs1\), independent of
\(\chi,p,n,\delta\), with the following property.
Put
\begin{equation}\label{eq:branching}
 K=C_3(r)\Lambda^{2p/r},
\end{equation}
and let \(t=t_{j-1}\) for some \(1\leqs j\leqs m_k\).
For every pair \((\sigma,\varphi)\) that is \(C^r\)-bounded
at time \(t\), with \(\sigma\) satisfying
\eqref{eq:original-coordinates}, there is a target chart
\(\varphi^+\) with the following
covering bounds for points of \(\F_c\) covered by \(\sigma\):
\begin{enumerate}
\item at most \(K\) maps for points with \(f^{\tau+t+p}x\in P\),
      for any prescribed atom \(P\in\PP_\delta\);
\item at most \(K\) maps for points with \(\wh z_t(x)\in\Gamma_p\);
\item at most \(K\Lambda^{2p}\) maps for all points covered by \(\sigma\).
\end{enumerate}
Each covering map has the form \(\sigma\circ\theta\), where
\(\theta:[0,1]^2\to[0,1]^2\) is smooth,
\(\|D\theta\|_\infty\leqs1\), and
\((\sigma\circ\theta,\varphi^+)\) is \(C^r\)-bounded at time \(t+p\).
\end{lemma}

\begin{proof}
Write
\(\psi=\varphi\circ V_t\circ\sigma\).
The bound \(|\psi|_r\leqs\delta\), the \(C_2\)-Lipschitz
bound for \(\varphi^{-1}\), and
\(\|Df^p\|_\infty\leqs\Lambda^p\) bound
the diameter of \(f^p(V_t(\sigma([0,1]^2)))\) by
\(\sqrt2 C_2\Lambda^p\delta\).
The scale condition therefore gives a chart \(\varphi^+\) such that
\[
 f^p(V_t(\sigma([0,1]^2)))
 \subset(\varphi^+)^{-1}((-2,2)^2).
\]
This chart is fixed before choosing the reparametrizations.
Set \(F_p=\varphi^+\circ f^p\circ\varphi^{-1}\) and define
the rescaled map
\[
 g=\delta^{-1}\varphi^+\circ V_{t+p}\circ\sigma
  =\delta^{-1}F_p\circ\psi.
\]
Since
\[
 |\varphi^+\circ V_{t+p}\circ\sigma\circ\theta|_r
 =\delta|g\circ\theta|_r,
\]
it suffices to cover the parameters corresponding
to the points in each case by maps \(\theta\) with
\(\|D\theta\|_\infty\leqs1\) and
\(|g\circ\theta|_r\leqs1\).
In the chain rule for \(F_p\circ\psi\), terms containing
\(DF_p\) have one derivative factor of \(\psi\), whereas
those containing \(D^jF_p\), \(j\geqs2\), have \(j\) such factors.
Since \(|\psi|_r\leqs\delta\leqs1\), division by \(\delta\)
bounds these terms by \(\|DF_p\|_\infty\) and
\(\delta|F_p|_r\), respectively.
Thus \eqref{eq:jets} gives
\begin{equation}\label{eq:block-jet}
 |g|_r
 \leqs c_r\max\{\|DF_p\|_\infty,\delta|F_p|_r\}
 \leqs c_r\Lambda^p,
\end{equation}
where \(c_r\geqs1\) depends only on \(r\).

For part~(1), fix \(P\in\PP_\delta\). The values \(g(v)\) at
parameters with \(\sigma(v)=u_x\) and \(f^{\tau+t+p}x\in P\)
have diameter at most
\(C_2\delta^{-1}\operatorname{diam}P\leqs C_2\).
If this set is nonempty, any one of these values is the
center of a ball of radius \(C_2\) containing all of them.
The value case \(S_2\) of Lemma~\ref{lem:reparametrization}
therefore applies.

For part~(2), take a covered point \(x\) with
\(\wh z_t(x)\in\Gamma_p\) and a parameter \(v\) with
\(\sigma(v)=u_x\).
Put \(A=D_{f^{\tau+t}x}f^p\), so that
\(AL_t(x)=L_{t+p}(x)\).
By \eqref{eq:good-projective-set} and the choice of \(\Lambda\),
\[
 \|A|_{L_t(x)}\|<1,\qquad \|A^{-1}\|\leqs\Lambda^p.
\]
Set
\[
 T=D_v(V_t\circ\sigma),\qquad
 N=\pi_{L_{t+p}(x)^\perp}AT.
\]
Then \(AT=D_v(V_{t+p}\circ\sigma)\).
The \(C^r\) bound at time \(t\) and
\eqref{eq:projection} at time \(t+p\), applied to the same
map \(\sigma\), give
\[
 \|T\|\leqs C_2\delta,\qquad
 \|N\|\leqs\delta\Lambda^{-p}.
\]
The image of \(AT-N\) lies in \(L_{t+p}(x)=AL_t(x)\),
so the image of \(T-A^{-1}N\) lies in \(L_t(x)\).
Using \(AT=A(T-A^{-1}N)+N\), we obtain
\begin{equation}\label{eq:normal-to-full}
 \begin{aligned}
 \|D_v(V_{t+p}\circ\sigma)\|
 &=\|AT\|\\
 &\leqs\|A|_{L_t(x)}\|
       \bigl(\|T\|+\|A^{-1}\|\,\|N\|\bigr)+\|N\|\\
 &\leqs C_2\delta+\Lambda^p\delta\Lambda^{-p}
       +\delta\Lambda^{-p}\leqs(C_2+2)\delta.
 \end{aligned}
\end{equation}
At each such parameter \(v\), passing to the target chart
and dividing by \(\delta\) gives
\[
 \|D_vg\|\leqs C_2(C_2+2).
\]
This is the first-derivative bound required by case \(S_1\)
of Lemma~\ref{lem:reparametrization}, and is independent of \(p\).

For both parts, apply Lemma~\ref{lem:reparametrization}
to \(g\), with \(B=C_2(C_2+2)\geqs C_2\) and
highest-derivative bound \(c_r\Lambda^p\).
This gives at most
\(C_1(r,B)c_r^{2/r}\Lambda^{2p/r}\)
reparametrizations. For part~(3), subdivide
\([0,1]^2\) into \(m^2\) equal squares, where
\(m=\lceil c_r\Lambda^p\rceil\).
Their affine parametrizations satisfy
\[
 \|D\theta\|_\infty=m^{-1}\leqs1,\qquad
 |g\circ\theta|_r\leqs m^{-1}|g|_r\leqs1,
 \qquad m^2\leqs4c_r^2\Lambda^{2p}.
\]
The three covering bounds follow by taking \(C_3(r)\) large
enough to dominate \(C_1(r,B)c_r^{2/r}\) and \(4c_r^2\).
In every case,
\(|\varphi^+\circ V_{t+p}\circ\sigma\circ\theta|_r\leqs\delta\), and
\(\|D\theta\|_\infty\leqs1\) preserves
\eqref{eq:original-coordinates}.
If \(\sigma(v)=u_x\) for a point \(x\) in the set being
covered, then \(v=\theta_i(v')\) for some \(i\) and
\(v'\in[0,1]^2\), so \((\sigma\circ\theta_i)(v')=u_x\).
\end{proof}

\subsection{Counting completions}

The remaining multiplicity will be bounded by the number of
periodic points in a Bowen ball. Set
\[
 M_n=\max\left\{1,\sup_{y\in M}
       \#\bigl(B_n(y,\sqrt2 C_2\Lambda^p\delta)
                    \cap\Per_\chi(n,M)\bigr)\right\}.
\]

\begin{lemma}[Completion bound at a fixed scale]
\label{lem:finite-scale-completion}
With the parameters above, for every \(n>p\), every grid
offset \(k\), and every erased-block selection \(\mathcal I(x)\),
\begin{align}
 \log\#\F_n-H_{\nu_{\F_n}}(C)\leqs{}&
 \left(\frac{2\log\Lambda}{r}
       +\frac{\log(2C_3(r))}{p}\right)n\notag\\*
 &+2p\log\Lambda\int N_n^{\mathrm{out}}\,d\nu_{\F_n}
 +\log M_n+\log(9C_2^2\Lambda^{2p}\delta^{-2}).
 \label{eq:completion-finite-scale}
\end{align}
\end{lemma}
\begin{proof}
Put \(\nu_n=\nu_{\F_n}\). We construct a sequence of
partitions \(\mathcal Q_0,\ldots,\mathcal Q_{m_k}\), each
refining its predecessor, with \(\mathcal Q_0\) refining the
partition into fibers of \(C\).
To each atom of \(\mathcal Q_j\), we will
assign a pair \((\sigma,\varphi)\) that is \(C^r\)-bounded
at time \(t_j\) and whose parameter map covers every point
of the atom.
For each covering family, fix an ordering and assign each
point to the first map that covers it.

\emph{(1) The initial partition.}
On each fiber \(\F_c\), write \(\Phi_0=\Phi_i\) and subdivide
\([-1,1]^2\) into \(m^2\) equal closed squares, where
\(m=\lceil2C_2\Lambda^p/\delta\rceil\).
Every affine parametrization \(\sigma\) of one of these squares satisfies
\[
 \varphi_i\circ V_0\circ\sigma=\sigma,\qquad
 s:=|\sigma|_r=\|D\sigma\|_\infty
 =2/m\leqs\frac{\delta}{C_2\Lambda^p}\leqs1.
\]
Thus \(\sigma\) satisfies \eqref{eq:original-coordinates}, and
\((\sigma,\varphi_i)\) is \(C^r\)-bounded at time zero.
If \(k>0\), then
\[
 \operatorname{diam}V_k(\sigma([0,1]^2))
 \leqs\sqrt2 C_2\Lambda^k s\leqs\sqrt2\delta.
\]
The scale condition gives a chart \(\varphi^+\) containing
this image in its inner domain.
Put \(F=\varphi^+\circ f^k\circ\Phi_0\).
Since \(\sigma\) is affine and \(s\leqs1\), the
first-derivative bound and \eqref{eq:jets} give
\[
 |\varphi^+\circ V_k\circ\sigma|_r
 =|F\circ\sigma|_r
 \leqs\max\{\Lambda^p s,\delta^{-1}s^2\}\leqs\delta.
\]
When \(k=0\), use \(\varphi^+=\varphi_i\).
Assigning points to these maps on each fiber defines
\(\mathcal Q_0\), with assigned pairs at time \(t_0=k\), and
\[
 H_{\nu_n}(\mathcal Q_0\mid C)
 \leqs\log m^2\leqs\log(9C_2^2\Lambda^{2p}\delta^{-2}).
\]

\emph{(2) Refinement along the grid.}
For \(1\leqs j\leqs m_k\), put
\[
 B_j=\{x\in\F_n:j\in\mathcal I(x),\
          \wh z_{t_{j-1}}(x)\notin\Gamma_p\}.
\]
Suppose \(\mathcal Q_{j-1}\) has been constructed and let
\(A\) be one of its atoms. Since \(A\) lies in a fiber of \(C\),
the block \(\{t_{j-1}+1,\ldots,t_j\}\) is either retained for
every point of \(A\) or erased for every point of \(A\).
Apply Lemma~\ref{lem:block-covering} to the pair assigned to \(A\),
fixing the target chart and the required covering families.
If the block is retained, use part~(1) with the endpoint atom
specified by \(C\). If it is erased, split \(A\) into
\(A\setminus B_j\) and \(A\cap B_j\), and use parts~(2) and~(3),
respectively.
Assign points within each group to the corresponding covering
maps, keeping the two groups separate even if their covering
images overlap. This defines \(\mathcal Q_j\) as a refinement
of \(\mathcal Q_{j-1}\), with \(C^r\)-bounded pairs at time \(t_j\).
Thus \(A\setminus B_j\) is divided into at most \(K\) atoms
of \(\mathcal Q_j\), and \(A\cap B_j\) into at most
\(K\Lambda^{2p}\).

Condition on \(\mathcal Q_{j-1}\), and then on
\(\one_{B_j}\). The binary entropy is at most \(\log2\);
the two covering bounds therefore give
\[
 \begin{aligned}
 H_{\nu_n}(\mathcal Q_j\mid\mathcal Q_{j-1})
 &\leqs\log2+(1-\nu_n(B_j))\log K
       +\nu_n(B_j)\log(K\Lambda^{2p})\\
 &=\log(2K)+2p\log\Lambda\,\nu_n(B_j).
 \end{aligned}
\]
By definition,
\(\sum_{j=1}^{m_k}\one_{B_j}=N_n^{\mathrm{out}}\).

\emph{(3) Final multiplicity and the counting bound.}
Fix \(A\in\mathcal Q_{m_k}\). The labels \(\tau,\Phi_0\) are
constant on \(A\). Each refinement preserves all preceding
assignments, so all points of \(A\) were assigned to the same
map at each sampled time \(s\in\mathcal T_k\), including time zero.
The \(C^r\) bounds at the sampled times and the
\(C_2\)-Lipschitz bound for the inverse charts therefore give
\[
 \operatorname{diam}f^{\tau+s}(A)\leqs\sqrt2 C_2\delta
 \qquad(s\in\mathcal T_k).
\]
For every \(0\leqs t<n\), the largest sampled time \(s\leqs t\)
satisfies \(0\leqs t-s\leqs p-1\), including after the last
complete block. Hence
\[
 \operatorname{diam}f^{\tau+t}(A)
 \leqs\sqrt2 C_2\Lambda^{p-1}\delta
 <\sqrt2 C_2\Lambda^p\delta.
\]
Choosing \(x_A\in A\), we obtain
\[
 f^\tau(A)\subset
 B_n(f^\tau x_A,\sqrt2 C_2\Lambda^p\delta)
 \cap\Per_\chi(n,M).
\]
Since \(f^\tau\) is injective, \(\#A\leqs M_n\).

Since \(\nu_n\) is uniform on \(\F_n\), its conditional
measure on each partition atom is uniform.
As \(\mathcal Q_{m_k}\) refines the partition into fibers of
\(C\), the counting identity \eqref{eq:name-decomposition}
and the entropy chain rule give
\[
 \begin{aligned}
 \log\#\F_n-H_{\nu_n}(C)
 &=H_{\nu_n}(\mathcal Q_{m_k}\mid C)
   +\sum_{A\in\mathcal Q_{m_k}}\nu_n(A)\log\#A\\
 &\leqs H_{\nu_n}(\mathcal Q_0\mid C)
   +\sum_{j=1}^{m_k}H_{\nu_n}(\mathcal Q_j\mid\mathcal Q_{j-1})
   +\log M_n\\
 &\leqs\log(9C_2^2\Lambda^{2p}\delta^{-2})+m_k\log(2K)\\
 &\quad+2p\log\Lambda\int N_n^{\mathrm{out}}\,d\nu_n
   +\log M_n.
 \end{aligned}
\]
Using \(m_k\leqs n/p\) and
\(K=C_3(r)\Lambda^{2p/r}\) yields
\eqref{eq:completion-finite-scale}.
\end{proof}

\begin{proof}[Proof of Proposition~\ref{prop:completion}]
Fix \(\varepsilon>0\).
Set \(C_4=2\log\Lambda\), independently of
\(\chi,\varepsilon,r,p\).
We apply Lemma~\ref{lem:finite-scale-completion}, choosing
\(r,p,\delta\) in that order.
Smoothness allows us to take \(r\) arbitrarily large; choose \(r\)
so that \(2\log\Lambda/r<\varepsilon/3\).
Since \(C_3(r)\) is independent of \(p\), choose
\(p_1=p_1(\varepsilon)\) so that every \(p\geqs p_1\)
satisfies
\[
 \frac{\log(2C_3(r))}{p}<\frac{\varepsilon}{3}.
\]
For any such fixed \(p\), local periodic growth
\eqref{eq:local} allows us to choose a fixed
\(0<\delta\leqs\delta_0(r,p)\) such that
\begin{equation}\label{eq:multiplicity}
 \limsup_{n\to\infty}\frac1n\log M_n\leqs\varepsilon/3.
\end{equation}
With these parameters fixed,
\(\log(9C_2^2\Lambda^{2p}\delta^{-2})\) is independent of \(n\).
The two strict bounds above and \eqref{eq:multiplicity}
therefore give \eqref{eq:completion-entropy} for all
sufficiently large \(n\).
The threshold is uniform over all grid offsets and
erased-block choices.
\end{proof}
\section{Global periodic counts}\label{sec:global}

We first extend the maximal-growth argument to finite mixtures,
adapting the retained occupation estimate and the concentration
of individual orbit measures.
Common Pesin blocks then give an exponentially smaller remainder
outside the periodic orbits homoclinically related to maximal
measures, allowing us to sum over disjoint counting sets.
For a saddle periodic orbit \(\mathcal O\) and \(x\in\mathcal O\),
write \(\nu_{\mathcal O}=\nu_x\) and
\(\wh\nu_{\mathcal O}=\wh\nu_x\) for its orbit measure and
unstable lift.

\subsection{Finite mixtures at maximal growth}

\begin{lemma}[Finite mixtures at maximal growth]
\label{lem:finite-mixture-growth}
Let \(\chi>0\), and let
\(\F_n\subset\Per_\chi(n,M)\) be nonempty finite invariant sets.
Suppose that
\[
 \nu_{\F_n}\longrightarrow\nu=\sum_i t_i\mu_i,
 \qquad \frac1n\log\#\F_n\longrightarrow h_0,
\]
where \(\nu\) is a finite convex combination of ergodic
invariant measures \(\mu_i\) of common entropy \(h_0>0\).
Put \(\mathcal S=\operatorname{conv}\{\wh\mu_i^u\}\).
For every \(\varepsilon>0\),
\[
 \frac{\#\{x\in\F_n:d_{\mathrm w}(\wh\nu_x,\mathcal S)
                         \geqs\varepsilon\}}{\#\F_n}
 \longrightarrow0.
\]
\end{lemma}

\begin{proof}
Pass to a subsequence on which the unstable lifts converge.
The classification of invariant lifts in Section~\ref{sec:prelim},
applied on the disjoint invariant sets of generic points of
the \(\mu_i\), gives a limit of the form
\[
 \sum_i\bigl((t_i-b_i)\wh\mu_i^u+b_i\wh\mu_i^s\bigr),
 \qquad 0\leqs b_i\leqs t_i.
\]
The proof of Proposition~\ref{prop:counting} extends to this
limit as follows. In Lemma~\ref{lem:retained}, choose the
partition with \(\nu(\partial\PP)=0\).
The retained occupation limit is invariant and bounded by
\(\nu\), so it has the form \(\eta=\sum_i s_i\mu_i\),
with \(0\leqs s_i\leqs t_i\).
Put \(t=\eta(M)\). If \(t>0\), affinity of entropy gives
\(h_{\eta/t}(f)=h_0\), and \(\eta\leqs\nu\) gives
\((\eta/t)(\partial\PP)=0\).
The window-entropy calculation in Lemma~\ref{lem:retained}
therefore still gives \eqref{eq:retained}; when \(t=0\),
the limiting window law is \(\delta_{\ast^m}\).
Equation~\eqref{eq:stable-block-mass} holds for every \(\mu_i\),
so \(p\) can be chosen simultaneously for this finite family.
The normal control from Subsection~\ref{subsec:normal} and
the completion estimates from Section~\ref{sec:geometry}
apply unchanged.
In the last step of the proof of Proposition~\ref{prop:counting},
replace \(\omega\wh\mu^s\) by \(\sum_i b_i\wh\mu_i^s\).
The latter measure is invariant, has mass \(\sum_i b_i\),
and assigns mass tending to zero to \(\Gamma_p^c\).
Consequently,
\[
 \limsup_{n\to\infty}\frac1n\log\#\F_n
 \leqs\left(1-\sum_i b_i\right)h_0.
\]
The growth assumption forces every \(b_i\) to vanish.
As every limit of the lifts equals
\(\wh\nu=\sum_i t_i\wh\mu_i^u\), we have
\(\wh\nu_{\F_n}\to\wh\nu\).

For the orbit-measure assertion, adapt the distribution argument
in the proof of Lemma~\ref{lem:orbit-measures} to
\((\wh M,\wh f)\).
Every limiting distribution \(\Theta\) of the lifted orbit
measures has barycenter \(\wh\nu\).
Let \(\kappa_\rho\) denote the ergodic decomposition of an
invariant measure \(\rho\). Uniqueness gives
\[
 \int\kappa_\rho\,d\Theta(\rho)
 =\sum_i t_i\delta_{\wh\mu_i^u}.
\]
Thus \(\kappa_\rho\) is supported on
\(\{\wh\mu_i^u\}\) for \(\Theta\)-almost every \(\rho\), so
\(\Theta(\mathcal S)=1\).
Since \(\mathcal S\) is compact, the Portmanteau theorem applied to
\(\{\rho:d_{\mathrm w}(\rho,\mathcal S)\geqs\varepsilon\}\)
proves the assertion.
\end{proof}

\subsection{Weighted Pesin blocks and homoclinic localization}

We use the homoclinic relation of \cite[\S2.4]{BCSMME}:
\(\mathcal O\sim_h\mu\), for an ergodic hyperbolic measure
\(\mu\) of saddle type, means that
\(W^u(\mathcal O)\) meets \(W^s(x)\) and
\(W^s(\mathcal O)\) meets \(W^u(x)\), both transversely,
for \(x\) in a positive \(\mu\)-measure set of regular points.
For periodic measures, this agrees with the relation defined
in the introduction.

The following finite family version of the localization
argument in \cite[Proposition~4.8]{LT} uses the weighted
Pesin blocks of~\cite{LTTY}, see also \cite{BCSSPR}.

\begin{lemma}[Homoclinic localization]
\label{lem:homoclinic-localization}
Let \(\mu_1,\ldots,\mu_N\) be ergodic hyperbolic measures
of saddle type, and let \(\mathcal O_j\) be saddle periodic
orbits. If
\[
 \wh\nu_{\mathcal O_j}\longrightarrow
       \sum_i t_i\wh\mu_i^u,
 \qquad t_i\geqs0,\qquad \sum_i t_i=1,
\]
then, for every sufficiently large \(j\),
\(\mathcal O_j\sim_h\mu_i\) for some \(i\).
\end{lemma}

\begin{proof}
Put \(\nu_j=\nu_{\mathcal O_j}\) and
\(\nu=\sum_i t_i\mu_i\). Projection gives
\(\nu_j\to\nu\). Integrating the continuous
functions
\[
 (x,L)\longmapsto\log\|D_xf|_L\|,
 \qquad x\longmapsto\log|\det D_xf|
\]
gives
\[
 \lambda^u(\nu_j)\longrightarrow
       \sum_i t_i\lambda^u(\mu_i),\qquad
 \lambda^s(\nu_j)\longrightarrow
       \sum_i t_i\lambda^s(\mu_i).
\]
For \(f^{-1}\), the positive exponent is \(-\lambda^s\).
Using these limits and the common positive lower bound on
the positive exponents of the \(\mu_i\) in both time directions,
apply \cite[Proposition~2.11]{LTTY} with \(k=1\)
to the ergodic measures \(\nu_j\) for \(f\) and \(f^{-1}\),
with exceptional mass less than \(1/4\) in each direction.
Let \(\Lambda_j^u,\Lambda_j^s\) be the resulting weighted
blocks. In each direction, the weight function, block
parameters, and iterate are independent of \(j\).
At points of \(\mathcal O_j\) in both blocks,
the splittings of both blocks agree with the Oseledets
splitting, with the two bundles interchanged for \(f^{-1}\).
These common choices allow us to take the closure
\[
 \Lambda=\overline{\bigcup_{j\geqs j_0}
       (\mathcal O_j\cap\Lambda_j^u\cap\Lambda_j^s)}
\]
and retain both splittings and their angle bounds
\cite[Propositions~2.1--2.2]{LTTY}.
This compact set satisfies \(\nu_j(\Lambda)>1/2\)
for all large \(j\). Its stable and unstable plaques have
uniform positive size, uniformly transverse tangent spaces,
and \(C^1\)-continuous dependence on the base point
\cite[Proposition~2.15]{LTTY}.
By \cite[Lemma~2.10]{LTTY} and Birkhoff's theorem,
these plaques are contained in the corresponding Pesin
manifolds at periodic points and at \(\mu_i\)-almost every
point of \(\Lambda\).

Cover \(\Lambda\) by finitely many compact subsets
\(D_1,\ldots,D_s\) small enough that any two points in
one set have transversely intersecting stable and unstable
plaques in both directions.
If the conclusion failed, pass to a subsequence with
\(\mathcal O_j\not\sim_h\mu_i\) for every \(i\).
Along a further subsequence, some fixed \(D=D_a\) satisfies
\(\nu_j(D)\geqs1/(2s)\). Since \(D\) is closed,
\[
 \nu(D)\geqs\limsup_j\nu_j(D)
       \geqs1/(2s)>0.
\]
Thus \(\mu_i(D)>0\) for some \(i\).
Since \(\mathcal O_j\cap D\ne\varnothing\), the choice of
\(D\) gives both transverse intersections with
\(\mu_i\)-almost every point of \(D\).
Hence \(\mathcal O_j\sim_h\mu_i\), a contradiction.
\end{proof}

\subsection{An entropy gap outside}

\begin{proposition}[Homoclinic localization of periodic counts]
\label{prop:global-remainder}
Let \(K\subset M\) be compact and invariant, with
\(h=\htop(f|K)>0\). Suppose that its ergodic maximal
measures form a finite nonempty family
\(\mu_1,\ldots,\mu_N\). Put
\[
 \mathcal C_i=\bigcup_{\mathcal O\sim_h\mu_i}\mathcal O,
 \qquad \mathcal C=\bigcup_i\mathcal C_i,
\]
where the first union is over saddle periodic orbits in \(M\).
For every \(\chi>0\), there are \(\eta>0\) and
\(n_0\geqs1\) such that
\begin{equation}\label{eq:global-remainder}
 \#\bigl(\Per_\chi(n,K)\setminus\mathcal C\bigr)
 \leqs e^{(h-\eta)n}\qquad(n\geqs n_0).
\end{equation}
\end{proposition}

\begin{proof}
Write \(\mathcal E_n=\Per_\chi(n,K)\setminus\mathcal C\).
Lemma~\ref{lem:periodic-entropy} bounds its upper growth rate
by \(h\). If equality held, we could pass to a sequence with
\[
 \frac1n\log\#\mathcal E_n\longrightarrow h,
 \qquad \nu_{\mathcal E_n}\longrightarrow\sum_i t_i\mu_i.
\]
Indeed, the same lemma makes every counting-measure limit
maximal on \(K\).
Lemma~\ref{lem:finite-mixture-growth} gives points
\(x_n\in\mathcal E_n\) whose unstable orbit measures approach
\(\operatorname{conv}\{\wh\mu_1^u,\ldots,\wh\mu_N^u\}\).
After taking a subsequence, these orbit measures converge to
a convex combination of these lifts.
Lemma~\ref{lem:homoclinic-localization} therefore places the
orbit of \(x_n\) in some \(\mathcal C_i\) for all large \(n\),
a contradiction. The upper growth rate is strictly below
\(h\), which proves \eqref{eq:global-remainder}.
\end{proof}

\subsection{Transfer from classes to the whole surface}

\begin{proof}[Proof of Corollary~\ref{cor:global}]
Fix \(\chi>0\) and write \(h=\htop(f)>0\).
By \cite[Theorems~1--2]{BCSMME}, the ergodic maximal
measures \(\mu_1,\ldots,\mu_N\) are finite in number and
pairwise non-homoclinically related. Their supports
\(\HH_i\) are homoclinic classes of entropy \(h\), with
unique maximal measures \(\mu_i\) and homoclinic periods
\(\ell_i\).
Let \(\mathcal C_i\) be the periodic counting sets in
Proposition~\ref{prop:global-remainder}. They are pairwise
disjoint, \(\mathcal C_i\subset\HH_i\), and every orbit
in \(\mathcal C_i\) has period divisible by \(\ell_i\).
These are the periodic parts of the Borel homoclinic classes
used in the irreducible codings of
\cite[\S\S2.4 and~3]{BCSMME} and~\cite[\S4]{LT}.

Apply Proposition~\ref{prop:global-remainder} to \(M\) and
to each \(\HH_i\). For some common \(\eta>0\),
\[
 \begin{aligned}
 \#\Per_\chi(n,M)
 &=\sum_i\#\bigl(\Per_\chi(n,M)\cap\mathcal C_i\bigr)
        +O(e^{(h-\eta)n}),\\
 \#\Per_\chi(n,\HH_i)
 &=\#\bigl(\Per_\chi(n,M)\cap\mathcal C_i\bigr)
        +O(e^{(h-\eta)n}).
 \end{aligned}
\]
The \(i\)-th summand in the first line vanishes unless
\(\ell_i\mid n\). Consequently,
\begin{equation}\label{eq:global-transfer}
 \#\Per_\chi(n,M)
 =\sum_{\ell_i\mid n}\#\Per_\chi(n,\HH_i)
       +O(e^{(h-\eta)n}).
\end{equation}
For every fixed \(\chi>0\), formula~\eqref{eq:global-transfer}
also gives, for some \(\eta>0\),
\begin{equation}\label{eq:global-off-period}
 \#\Per_\chi(n,M)=O(e^{(h-\eta)n})
 \qquad\text{if }\ell_i\nmid n\text{ for every }i.
\end{equation}
For \(0<\chi<h\), apply Theorem~\ref{thm:prime} to each
term on the right to obtain
\[
 e^{-nh}\#\Per_\chi(n,M)
       =\sum_{\ell_i\mid n}\ell_i+o(1).
\]
Topological transitivity gives \(N=1\), and topological
mixing also gives \(\ell_1=1\), by
\cite[Theorems~1--2]{BCSMME}.
\end{proof}

\section{Dynamical zeta functions}\label{sec:zeta}

The class-to-global comparison \eqref{eq:global-transfer}
already has an exponential error. We now obtain the same
precision relative to the explicit main term within each
class, then show that sinks and sources have growth strictly
below the full entropy. Together these estimates prove
Theorem~\ref{thm:zeta}.

\subsection{Exponential counting on a homoclinic class}

\begin{proposition}[Exponential remainder on a class]
\label{prop:class-zeta-remainder}
Let \(\HH\) be a closed homoclinic class of entropy \(h>0\)
and homoclinic period \(\ell\). For every \(0<\chi<h\),
there is \(\eta>0\) such that
\[
 \#\Per_\chi(n,\HH)
 =\ell\one_{\{\ell\mid n\}}e^{hn}+O(e^{(h-\eta)n}).
\]
\end{proposition}

\begin{proof}
\emph{Comparison with cylinder counts.}
Let \(\mu\) be the unique maximal measure on \(\HH\).
By \eqref{eq:cutoff-comparison}, it suffices to prove the
formula for the reference cutoffs \(\alpha_*,\beta_*\)
from \eqref{eq:LT}.

The proof of \cite[Proposition~4.5]{LT} applies without the
mixing reduction and gives an irreducible SPR shift \((Z,\sigma)\)
of entropy \(h\), together with a semiconjugacy
\(\Pi:Z\to\HH\) mapping its maximal measure
isomorphically to \(\mu\). The finite extension in that
construction gives a common hyperbolicity cutoff
\(\kappa>0\) for all periodic images.
The Bernoulli-times-cycle descriptions in
\cite[Theorem~6.6(1)]{BCSSPR} and
\cite[Theorem~2(1)]{BCSMME} imply that \(Z\) has period
\(\ell\).

For a nonempty cylinder \(C\subset Z\), let \(S_n(C)\)
count the points of least period \(n\) whose orbits visit
\(C\). We first find nonempty cylinders \(C_-,C_+\subset Z\)
and \(\eta>0\) such that
\begin{equation}\label{eq:class-cylinder-comparison}
 S_n(C_-)-O(e^{(h-\eta)n})
 \leqs\#\Per_{\alpha_*,\beta_*}(n,\HH)
 \leqs S_n(C_+)+O(e^{(h-\eta)n}).
\end{equation}
For the upper bound, the proofs of
\cite[Propositions~2.4 and~4.10, Corollary~5.2]{LT}
apply to the original coding and its Bowen quotient and give a cylinder
\(C_+\) such that all
but \(O(e^{(h-\eta)n})\) points in
\(\Per_{\alpha_*,\beta_*}(n,\HH)\) have an
\(n\)-periodic \(\Pi\)-preimage whose orbit visits
\(C_+\). Each such preimage has least period \(n\),
since its image does.

For the lower bound, the same construction of \(Z\) gives
a finite-to-one coding on its regular part, with a locally
finite Bowen relation of recurrent degree one; see the proof of
\cite[Proposition~4.5]{LT}. Hence
\cite[Theorem~4.6]{Buzzi} supplies a magic cylinder
\(C_-\) such that \(\Pi\) is injective on the periodic
points whose orbits visit it: every word along a periodic sequence recurs in
both time directions, and its Bowen fiber is finite.
This also preserves least periods.
Since \(\alpha_*<\lambda^u(\mu)\) and
\(\beta_*<-\lambda^s(\mu)\),
Theorem~\ref{thm:concentration}, applied with cutoff
\(\kappa\), bounds the number of images failing either
reference cutoff by \(O(e^{(h-\eta)n})\), after decreasing
\(\eta\). Injectivity proves the lower bound in
\eqref{eq:class-cylinder-comparison}.

\emph{Cylinder counts from first returns.}
We use the first return calculation of
\cite[Proposition~4.3]{BBG} and \cite[Proposition~9.2]{GS}.
Fix either cylinder \(C\in\{C_-,C_+\}\), and pass to a
higher-block presentation in which \(C\) is a vertex
cylinder. Let \(a_n(C)\) be the number of first return
loops of length \(n\) at this vertex, and put
\[
 F_C(z)=\sum_{n\geqs1}a_n(C)z^n.
\]
SPR gives a radius of convergence greater than \(e^{-h}\),
and positive recurrence gives \(F_C(e^{-h})=1\);
see \cite[\S2]{BBG} and \cite[proof of Lemma~5.3]{LT}.
The return lengths have gcd \(\ell\). Nonnegativity of
the coefficients therefore shows that \(1-F_C\) has no
zeros in \(|z|<e^{-h}\), and its zeros on
\(|z|=e^{-h}\) are precisely
\(e^{-h}\omega\), where \(\omega^\ell=1\).
These zeros are simple, by the \(\ell\)-fold rotational
symmetry and \(F_C'(e^{-h})>0\). Consequently, for some
\(R_C>e^{-h}\),
\[
 1-F_C(z)=(1-e^{h\ell}z^\ell)G_C(z),
 \qquad |z|<R_C,
\]
where \(G_C\) is holomorphic and nowhere zero.

Let \(Q_n(C)\) count the points fixed by \(\sigma^n\)
whose orbits visit \(C\). The first return decomposition
gives the following identity of formal power series,
as in \cite[proof of Lemma~5.3]{LT}:
\[
 \sum_{n\geqs1}Q_n(C)z^n
 =\frac{zF_C'(z)}{1-F_C(z)}
 =\frac{\ell e^{h\ell}z^\ell}{1-e^{h\ell}z^\ell}
       -z\frac{G_C'(z)}{G_C(z)}.
\]
Cauchy's coefficient estimate for the holomorphic term
\(zG_C'/G_C\) yields, for some \(\eta_C>0\),
\[
 Q_n(C)=\ell\one_{\{\ell\mid n\}}e^{hn}
          +O(e^{(h-\eta_C)n}).
\]
Since \(Q_n(C)=\sum_{d\mid n}S_d(C)\) and every proper
divisor of \(n\) is at most \(n/2\),
\[
 0\leqs Q_n(C)-S_n(C)
 \leqs\sum_{d\leqs n/2}Q_d(C)
 =O(e^{hn/2}).
\]
Decreasing \(\eta>0\) to work for both cylinders gives
\begin{equation}\label{eq:cylinder-exponential-count}
 S_n(C_\pm)
 =\ell\one_{\{\ell\mid n\}}e^{hn}
       +O(e^{(h-\eta)n}).
\end{equation}
Substitution into \eqref{eq:class-cylinder-comparison}
proves the formula for the reference cutoffs.
\end{proof}
\subsection{Sinks and sources}
\label{subsec:ss}

For \(\chi>0\), put
\[
 \begin{aligned}
 S^-_\chi(n)&=\{x\in\Fix(f^n):\rho(D_xf^n)\leqs e^{-\chi n}\},\\
 S^+_\chi(n)&=\{x\in\Fix(f^n):\rho(D_xf^{-n})\leqs e^{-\chi n}\},
 \end{aligned}
\]
where \(\rho\) denotes the spectral radius. For \(x\in\Fix(f^n)\),
the maximal Lyapunov exponent of \(f\) is
\(n^{-1}\log\rho(D_xf^n)\). These invariant sets
are finite: each is compact and consists of isolated fixed
points of \(f^n\).

\begin{proposition}[An entropy gap for sinks and sources]
\label{prop:ss-gap}
Suppose \(h=\htop(f)>0\). For every \(\chi>0\), there
are \(C,\eta>0\) such that
\begin{equation}\label{eq:ss-global-gap}
 \#S^-_\chi(n)+\#S^+_\chi(n)\leqs Ce^{(h-\eta)n}
 \qquad(n\geqs1).
\end{equation}
\end{proposition}

We prove the sink estimate; time reversal gives the
source estimate. The erasure argument of
Sections~\ref{sec:proofs} and~\ref{sec:geometry} turns a
positive limiting mass of contracting directions into a
strict entropy loss. We first prove the local periodic growth estimate
for sinks needed for that argument.

\begin{lemma}[Local growth of sinks]
\label{lem:ss-local-growth}
For every \(\chi>0\),
\begin{equation}\label{eq:ss-local-growth}
 \lim_{\epsilon\downarrow0}\limsup_{n\to\infty}\frac1n
 \log\max\left\{1,\sup_{x\in M}
 \#\bigl(S^-_\chi(n)\cap B_n(x,\epsilon)\bigr)\right\}=0.
\end{equation}
\end{lemma}

\begin{proof}
We group periodic points according to the growth of two directions
along successive orbit blocks. Reparametrization makes these growth
rates comparable throughout each image, transferring the return
contraction from a periodic point to the whole image. A subdivision
into pieces with short connecting paths then gives at most one fixed
point per piece.

\emph{(1) Reduction to norm contraction.}
Let \(\Lambda\) be as in Section~\ref{sec:geometry}, and put
\[
 T_a(n)=\{y\in\Fix(f^n):\|D_yf^n\|\leqs e^{-an}\},
 \qquad a>0.
\]
For every real \(2\times2\) matrix \(P\) and \(m\geqs2\),
complex Schur triangularization gives
\[
 \|P^m\|\leqs\rho(P)^m+m\rho(P)^{m-1}\|P\|.
\]
Indeed, in a unitary triangular form
\(\left(\begin{smallmatrix}\lambda_1&b\\0&\lambda_2\end{smallmatrix}\right)\)
of \(P\), the off-diagonal entry of the \(m\)-th power is
\(b\sum_{j=0}^{m-1}\lambda_1^{m-1-j}\lambda_2^j\),
and \(|b|\leqs\|P\|\).
Choose a fixed integer \(m\geqs2\) such that
\((m-1)\chi-\log\Lambda>m\chi/2\).
For \(y\in S^-_\chi(n)\), periodicity and
\(\|D_yf^n\|\leqs\Lambda^n\) give
\[
 \|D_yf^{mn}\|=\|(D_yf^n)^m\|
 \leqs e^{-m\chi n}
       +m e^{-((m-1)\chi-\log\Lambda)n}
 \leqs e^{-m\chi n/2}
\]
for all sufficiently large \(n\). Thus
\(S^-_\chi(n)\subset T_{\chi/2}(mn)\).

If \(S^-_\chi(n)\cap B_n(x,\epsilon)\) is nonempty, choose
\(z\) in this set. The \(n\)-step and \(mn\)-step Bowen
distances agree between \(n\)-periodic points, so
\[
 S^-_\chi(n)\cap B_n(x,\epsilon)
 \subset T_{\chi/2}(mn)\cap B_{mn}(z,2\epsilon).
\]
Since \(m\) is fixed, it suffices to prove zero local exponential
growth for \(T_a(n)\), for every \(a>0\), uniformly over centers
\(z\in\Fix(f^n)\).

\emph{(2) Reparametrizations with prescribed directional growth.}
Fix \(a>0\) and integers \(r\geqs2\), \(p>4/a\), and set
\[
 N=\lceil n/p\rceil,\qquad t_j=\min\{jp,n\}.
\]
For \(z\in\Fix(f^n)\), choose the charts of
Subsection~\ref{subsec:normal} with \(f^{t_j}z\in U'_{i_j}\), and put
\[
 \kappa_j(v)=\Phi_{i_j}\bigl(
   \varphi_{i_j}(f^{t_j}z)+C_2\epsilon v\bigr),
 \qquad 0\leqs j\leqs N,
\]
\[
 F_j=\kappa_{j+1}^{-1}\circ f^{t_{j+1}-t_j}\circ\kappa_j,
 \qquad 0\leqs j<N.
\]
Choose \(i_N=i_0\), so that \(\kappa_N=\kappa_0\).
Each block has length at most \(p\), and \(F_j(0)=0\).
For fixed \(r,p\), the relevant coordinate expressions of
\(f^k\), \(1\leqs k\leqs p\), have uniformly bounded derivatives
on the compact chart domains. Rescaling multiplies derivatives of
order \(q\geqs2\) by \((C_2\epsilon)^{q-1}\).
Hence, for \(0<\epsilon\leqs\epsilon(r,p)\) sufficiently small,
uniformly in \(n,z,j\), the maps \(F_j\) are defined on \(B(0,4)\) and
\[
 \kappa_j^{-1}\bigl(B(f^{t_j}z,\epsilon)\bigr)\subset B(0,1),
 \qquad
 \max_{2\leqs q\leqs r}\|D^qF_j\|_{\infty,B(0,4)}\leqs1.
\]
Multiply by a fixed cutoff supported in \(B(0,4)\) and equal to
one on \(B(0,3)\). Since \(F_j(0)=0\), this gives extensions to
\(\R^2\) satisfying
\[
 \|F_j\|_{C^r}\leqs C_r\Lambda^p,
 \qquad
 \|DF_j\|_\infty+
 \sup_{B(0,3)}\|(DF_j)^{-1}\|\leqs C\Lambda^p.
\]
Here \(C\) depends only on \(f\), the metric and the atlas,
and \(C_r\) may also depend on \(r\), but not on \(p\).
Write \(F^0=\mathrm{id}\) and \(F^j=F_{j-1}\circ\cdots\circ F_0\).

Fix an orthonormal basis \(e_1,e_2\).
For \(u\) with \(F^ju\in B(0,3)\) for \(0\leqs j<N\), and
\(v\in S^1\), define the one-block expansion
\[
 \ell_j(u,v)
 =\frac{\|D_uF^{j+1}v\|}{\|D_uF^jv\|}
 =\left\|D_{F^ju}F_j
     \frac{D_uF^jv}{\|D_uF^jv\|}\right\|.
\]
For \(y\in T_a(n)\cap B_n(z,\epsilon)\), let
\(u_y=\kappa_0^{-1}y\). Then
\[
 F^ju_y=\kappa_j^{-1}(f^{t_j}y)\in B(0,1)
 \quad(0\leqs j\leqs N),
 \qquad F^Nu_y=u_y.
\]
Record
\[
 k_{i,j}=\lfloor\log\ell_j(u_y,e_i)\rfloor,
 \qquad i=1,2,\quad 0\leqs j<N,
\]
so that \(e^{k_{i,j}}\leqs\ell_j(u_y,e_i)<e^{k_{i,j}+1}\).
The derivative bounds give
\[
 C^{-1}\Lambda^{-p}\leqs\ell_j(u_y,e_i)\leqs C\Lambda^p,
 \qquad |k_{i,j}|\leqs p\log\Lambda+C.
\]
Thus there are at most \((Cp)^{2N}\) pairs of records.

Fix a pair of records. The corresponding triples
\((u_y,e_1,e_2)\) satisfy the orbit and growth conditions in
\cite[Lemma~7.2]{Burguet}. Apply that lemma and
\cite[Remark~7.2]{Burguet} to \(F_0,\ldots,F_{N-1}\), with
the auxiliary integer \(s=r+1\) and a fixed cone parameter
giving aperture less than \(\pi/8\). In its notation, take
\[
 a_j=k_{1,j},\qquad b_j=-k_{2,j}-1,
\]
so that the second growth interval
\([e^{-b_j-1},e^{-b_j}]\) is \([e^{k_{2,j}},e^{k_{2,j}+1}]\).
We obtain a finite family of Nash maps, that is, real-analytic
semialgebraic maps, of the form
\[
 \psi(x,s_1,s_2)=
 \bigl(\psi_1(x),\psi_2(x,s_1),\psi_3(x,s_2)\bigr),
 \qquad (x,s_1,s_2)\in(0,1)^2\times(0,1)\times(0,1),
\]
continuously extendable to the closed parameter cube.
Their full images
\[
 \widehat E=\psi([0,1]^2\times[0,1]\times[0,1])
 \subset\R^2\times S^1\times S^1
\]
cover all triples \((u_y,e_1,e_2)\) with the chosen records.
Retain only full images meeting such a triple, and write
\(E=\psi_1([0,1]^2)\) for the corresponding base image.
The resulting family has the following properties:
\begin{enumerate}
\item[(i)] By the triangular construction in
  \cite[proof of Lemma~7.2 and Remarks~6.1, 7.2]{Burguet}, each \(\psi_1:(0,1)^2\to\R^2\) is a Nash
  diffeomorphism onto its image, and
  \[
   F^j(E)\subset B(0,2),\qquad 0\leqs j<N.
  \]
\item[(ii)] For every \((u,v_1,v_2)\in\widehat E\),
  \[
   e^{k_{i,j}-1}\leqs\ell_j(u,v_i)\leqs e^{k_{i,j}+2},
   \qquad i=1,2,\quad 0\leqs j<N.
  \]
  The two initial directions range in fixed cones of aperture
  less than \(\pi/8\). Since the full image meets a recorded
  triple, these cones contain \(e_1,e_2\), respectively; hence
  \[
   |\det(v_1,v_2)|\geqs\frac12
   \qquad\bigl((u,v_1,v_2)\in\widehat E\bigr).
  \]
\item[(iii)] For each fixed pair of records,
  \[
   \#\{\widehat E\}
   \leqs B_{r,p}A_r^N\Lambda^{8pN/(r-1)},
   \qquad \deg\psi_1\leqs A_r^N.
  \]
  Here \(A_r\geqs1\) is independent of \(p,n,z\), and
  \(B_{r,p}\geqs1\) is independent of \(n,z\); the bounds are
  uniform in the records and in \(0<\epsilon\leqs\epsilon(r,p)\).
\end{enumerate}
Here the degree of a scalar Nash function is the least total degree
of a nonzero polynomial vanishing on its graph; for a vector-valued
map, we take the maximum over its coordinates, as in
\cite[\S6.1]{Burguet}.

For the cardinality bound, the derivative bounds and the original records give
\[
 \max\left\{\frac{\|F_j\|_{C^r}}{e^{k_{1,j}}},
 \frac{\|F_j\|_{C^r}}{e^{k_{2,j}+1}}\right\}
 \leqs C_r'\Lambda^{2p}.
\]
The squared base factor in \cite[Lemma~7.2(iii)]{Burguet}
is therefore at most \(C_r\Lambda^{4p/(r-1)}\), including
its term \(\|F_j\|_{C^r}^{2/r}\); the two directional factors
together are at most \(C_r\Lambda^{4p/(r+1)}\).
Thus each block contributes at most
\[
 A_r\Lambda^{4p/(r-1)+4p/(r+1)}
 \leqs A_r\Lambda^{8p/(r-1)}.
\]
The constants absorbed into \(A_r\) depend on \(r\) but not
on \(p\). The prefactor in the cited lemma depends on the
uniform \(C^r\) bound \(C_r\Lambda^p\) and the fixed cone
parameter, and is absorbed into \(B_{r,p}\).

\emph{(3) Contraction on each image.}
Fix a retained full image \(\widehat E\), with the base image \(E\),
and put \(R=F^N\). Since all intermediate images of \(E\) lie
in \(B(0,2)\), 
\[
 R=\kappa_0^{-1}\circ f^n\circ\kappa_0
\]
on a neighborhood of \(E\).
Choose a recorded triple \((u_y,e_1,e_2)\in\widehat E\).
The chart bounds and \(y\in T_a(n)\) give
\(\|D_{u_y}R\|\leqs C_2^2e^{-an}\).
For any \(u\in E\), choose a lift
\((u,v_1,v_2)\in\widehat E\). Comparing with the original
records in this same full image yields
\[
 \begin{aligned}
 \|D_uR\,v_i\|
 &=\prod_{j=0}^{N-1}\ell_j(u,v_i)\\
 &\leqs\prod_{j=0}^{N-1}e^{k_{i,j}+2}
 \leqs e^{2N}\prod_{j=0}^{N-1}\ell_j(u_y,e_i)\\
 &=e^{2N}\|D_{u_y}R\,e_i\|
 \leqs C_2^2e^{-an+2N},\qquad i=1,2.
 \end{aligned}
\]
The determinant bound implies
\(\|D_uR\|\leqs C\max_i\|D_uR\,v_i\|\).
Thus the loss in transferring contraction is only \(Ce^{2N}\),
one fixed factor per block. Since \(N=\lceil n/p\rceil\) and
\(p>4/a\),
\[
 -an+2N\leqs-\frac a2n+2,
 \qquad
 \sup_E\|DR\|\leqs C C_2^2e^{-an+2N}<\frac1{\sqrt2}
\]
for all sufficiently large \(n\), uniformly over the retained images.

\emph{(4) Counting fixed points.}
To turn the derivative bound into a distance bound, we use the
geometric decomposition in the proof of
\cite[Lemma~6.2, Claims~1--2]{Burguet}, with the constant
horizontal and vertical fields. For the Nash diffeomorphism
\(\psi_1\) in part~(2), continuously extended to the closed square,
this gives a cover of \(E\) by at most
\[
 C(\deg\psi_1)^C\leqs CA_r^{CN}
\]
compact pieces \(Q\subset E\), where \(C\) is universal.
Any two points \(u,v\in Q\) can be joined within \(Q\) by a
rectifiable path \(\gamma\) with both coordinates monotone, and hence
\[
 \operatorname{length}(\gamma)
 \leqs |u_1-v_1|+|u_2-v_2|\leqs\sqrt2\,|u-v|.
\]

For distinct \(u,v\) in one such \(Q\), integration along a
connecting path \(\gamma\) now gives
\[
 \begin{aligned}
 |R(u)-R(v)|
 &\leqs\sup_E\|DR\|\operatorname{length}(\gamma)\\
 &\leqs\sqrt2\sup_E\|DR\|\,|u-v|<|u-v|.
 \end{aligned}
\]
Thus each piece contains at most one fixed point of \(R\), including
boundary points. Every \(y\in T_a(n)\cap B_n(z,\epsilon)\)
gives such a fixed point \(u_y\) in one of the retained base images.

Multiplying the number of records, the number of images per record,
and the number of pieces per image, we obtain
\[
 \begin{aligned}
 \#\bigl(T_a(n)\cap B_n(z,\epsilon)\bigr)
 &\leqs (Cp)^{2N}
       B_{r,p}A_r^N\Lambda^{8pN/(r-1)}\,C A_r^{CN}\\
 &\leqs C_{r,p}
       \left(A_rp^2\Lambda^{8p/(r-1)}\right)^{\lceil n/p\rceil},
 \end{aligned}
\]
after increasing \(A_r\), still independently of \(p,n,z\).
This bound holds for all sufficiently large \(n\), uniformly in \(z\).
The upper local exponential growth rate at this scale is therefore
at most
\[
 \frac{8\log\Lambda}{r-1}
 +\frac{\log A_r+2\log p}{p}.
\]
First choose \(r\) large, then \(p>4/a\) large, and finally
\(0<\epsilon\leqs\epsilon(r,p)\) small. The bound can be made
arbitrarily small. The fixed-power reduction in part~(1)
then proves \eqref{eq:ss-local-growth}.
\end{proof}

\begin{proof}[Proof of Proposition~\ref{prop:ss-gap}]
If the sink count had no exponential gap below \(h\),
we could choose a sequence of periods \(n\to\infty\) such that
\(\F_n=S^-_\chi(n)\) is nonempty and
\[
 \liminf_n\frac1n\log\#\F_n\geqs h.
\]

The return derivative of a sink may have nonreal eigenvalues.
Use the bundle \(\mathbb P_{\mathbb C}(TM\otimes_{\R}\mathbb C)\)
of complex lines, with the Hermitian extension of the metric,
and denote its projection and derivative action again by
\(\pi\) and \(\wh f\).
On each orbit in \(\F_n\), choose an eigenline for an
eigenvalue of largest modulus of its least-period return,
and transport it along the orbit. The resulting lines satisfy
\(D_xf\,L(x)=L(fx)\) and
\begin{equation}\label{eq:ss-line-quotient}
 \|D_xf^n|_{L(x)}\|=\rho(D_xf^n)\leqs e^{-\chi n},
 \qquad
 \frac{|\det D_xf^n|}{\|D_xf^n|_{L(x)}\|}
       \leqs e^{-\chi n}.
\end{equation}
The first inequality supplies the contracting mass, and
the second supplies the normal bound.

The counting lifts
\[
 \wh\nu_n=\frac1{\#\F_n}\sum_{x\in\F_n}\delta_{(x,L(x))}
\]
are invariant and project to \(\nu_{\F_n}\).
Pass to a subsequence with \(\wh\nu_n\to\wh\nu\), and
put \(\nu=\pi_*\wh\nu\).
For the continuous function \(\phi(x,L)=\log\|D_xf|_L\|\),
invariance gives
\begin{equation}\label{eq:ss-negative-mean}
 \int\phi\,d\wh\nu
 =\lim_n\frac1{n\#\F_n}
       \sum_{x\in\F_n}\log\rho(D_xf^n)
 \leqs-\chi.
\end{equation}
By Birkhoff's theorem, the restriction \(\wh\nu^-\) of
\(\wh\nu\) to the invariant set
\[
 \left\{(x,L):
   \lim_{p\to\infty}\frac1p\log\|D_xf^p|_L\|<0\right\}
\]
has mass \(\omega>0\).
Define \(\Gamma_p\) by \eqref{eq:good-projective-set} on
the complex line bundle. It is open, and bounded convergence
gives \(\wh\nu^-(\Gamma_p^c)\to0\).
We use this contracting mass in the erasure argument,
checking the geometric and entropy estimates below.

For the geometric estimate, the second inequality in
\eqref{eq:ss-line-quotient} says that \(J\circ\pi-\phi\)
has negative sum over each period.
The equality \eqref{eq:normal-area} holds in unitary bases,
so the proof of Lemma~\ref{lem:normal} gives a cyclic
origin \(\tau(x)\). Set
\(L_t(x)=L(f^{\tau(x)+t}x)\) in
\eqref{eq:transported-lifts}.
Then \eqref{eq:normal} holds with \(L_t(x)\) in place
of the unstable line, and \eqref{eq:transported-limit}
holds with \(\wh\nu_n\).
The coordinate labels are unchanged, and these constructions
require only \(f^nx=x\).

For completion, complexify the real operators \(A,T\)
in the proof of \eqref{eq:normal-to-full}, without changing
their operator norms.
With \(N=\pi_{L_{t+p}(x)^\perp}AT\) defined by Hermitian
projection, the image of \(T-A^{-1}N\) still lies in
\(L_t(x)\). The same calculation therefore gives
\eqref{eq:normal-to-full} for the original real operator
\(AT\). All parameter maps remain real and two-dimensional,
so the covering costs are unchanged.
Lemma~\ref{lem:ss-local-growth} supplies the final additional
local bound in Proposition~\ref{prop:completion}.

For the entropy estimate, choose the partition with
\(\nu\)-null boundary. In the contradiction argument of
Lemma~\ref{lem:retained}, every retained occupation limit
\(\eta\leqs\nu\) is invariant.
When \(\eta\ne0\), its normalization has entropy at most
\(h\) by the variational principle.
The window-entropy estimate thus gives \eqref{eq:retained}
with \(h_0=h\), without requiring \(\nu\) to be ergodic.

The proof of Lemma~\ref{lem:selection} now gives
\eqref{eq:selection-counting} with \(h_0=h\) and
\(\wh\nu_n\) in place of \(\wh\nu_{\F_n}\), using the
same definition of \(\Psi_{p,q}\) on the complex line bundle.
Since \(\wh\nu^-\) is invariant and
\(\wh\nu^-(\Gamma_p^c)\to0\), the final argument in the
proof of Proposition~\ref{prop:counting}, with
\(\omega\wh\mu^s\) replaced by \(\wh\nu^-\), yields
\[
 \limsup_n\frac1n\log\#\F_n\leqs(1-\omega)h<h,
\]
contrary to the choice of \(\F_n\).
Apply the same argument to \(f^{-1}\) for sources, and
increase \(C\) to cover the finitely many remaining periods.
\end{proof}

\subsection{Global counts and the zeta function}

\begin{proof}[Proof of Theorem~\ref{thm:zeta}]
Let \(h=\htop(f)>0\), fix \(0<\chi<h\), and write
\(\HH_i=\operatorname{supp}\mu_i\).
Apply Proposition~\ref{prop:class-zeta-remainder} to each
\(\HH_i\) and use the global transfer formula
\eqref{eq:global-transfer}. For some \(\eta_1>0\),
\begin{equation}\label{eq:saddle-exponential-remainder}
 \#\Per_\chi(n,M)
 =e^{hn}\sum_{i=1}^N\ell_i\one_{\{\ell_i\mid n\}}
       +O(e^{(h-\eta_1)n}).
\end{equation}
The \(\chi\)-hyperbolic saddle points fixed by \(f^n\) are counted by
\(\sum_{d\mid n}\#\Per_\chi(d,M)\).
Since every proper divisor of \(n\) is at most \(n/2\),
their contribution beyond least period \(n\) is
\[
 O\left(\sum_{d\leqs n/2}e^{hd}\right)=O(e^{hn/2}).
\]
Proposition~\ref{prop:ss-gap} bounds the remaining
\(\chi\)-hyperbolic sources and sinks by
\(O(e^{(h-\eta_2)n})\), for some \(\eta_2>0\).
Choosing \(0<\eta<\min\{\eta_1,\eta_2,h/2\}\) proves
\eqref{eq:hyperbolic-exponential-remainder}.

The error series
\[
 H_\chi(z)=\sum_{n\geqs1}\frac{z^n}{n}
 \left(\#\Fix_\chi(n,M)
       -e^{hn}\sum_i\ell_i\one_{\{\ell_i\mid n\}}\right)
\]
is holomorphic on \(|z|<e^{-h+\eta}\).
For \(|z|<e^{-h}\), absolute convergence and the power-series
expansion of \(-\log(1-w)\) give
\[
 \log\zeta_f^\chi(z)
 =H_\chi(z)-\sum_i\log(1-e^{h\ell_i}z^{\ell_i}),
\]
with the logarithms equal to zero at the origin.
Exponentiating proves \eqref{eq:hyperbolic-zeta-factorization}
and its meromorphic continuation.
The numerator never vanishes, and each denominator factor
has simple zeros. The stated pole orders follow, and in particular,
the pole at \(e^{-h}\) makes the Taylor radius exactly \(e^{-h}\).
\end{proof}

\end{document}